\documentclass[11pt,a4paper]{amsart}

\usepackage{amssymb,amscd}
\usepackage[all]{xy}
\allowdisplaybreaks

\newtheorem{theorem}{Theorem}[section]
\newtheorem{corollary}[theorem]{Corollary}
\newtheorem{proposition}[theorem]{Proposition}
\theoremstyle{definition}
\newtheorem{definition}[theorem]{Definition}
\newtheorem{example}[theorem]{Example}

\numberwithin{equation}{section}

\def\xarr{\xrightarrow}	\def\ti{\widetilde}
\def\CW{CW-complex} 	\def\CWs{CW-complexes}
\def\dlim{\varinjlim}	\def\+{\oplus}
\def\hot{\mathop\mathrm{Hot}\nolimits}
\def\hos{\mathop\mathrm{Hos}\nolimits}
\def\Hos{\mathbf{Hos}}	\def\xx{\times}
\def\wee{\wedge}		\def\sb{\subset}
		\def\b+{\bigoplus}
\def\set#1{\left\{\,#1\,\right\}}
\def\setsuch#1#2{\left\{\,#1\,|\,#2\,\right\}}
\def\iso{\simeq}		\def\cw{\mathsf{CW}}
\def\cwf{\mathsf{CWF}}	
\def\cws{\mathsf{CWS}}	\def\larr{\longrightarrow}
\def\iff{if and only if }		\def\str#1{\stackrel{#1}\larr}
\def\8{\infty}			\def\cel{*=0{\bullet}}
\def\im{\mathop\mathrm{Im}\nolimits}
\def\ker{\mathop\mathrm{Ker}\nolimits}
\def\GL{\mathop\mathrm{GL}\nolimits}
\def\mod{\,\mathop\mathrm{mod}}
\def\ds{\displaystyle}	\def\car{\ar@{-}}
\def\El{\mathop\mathbf{El}\nolimits}
\def\ob{\mathop\mathrm{ob}\nolimits}
\def\Ab{\mathbf{Ab}}	\def\*{\otimes}
\def\lst#1#2{ #1_1 , #1_2 , \dots , #1_{#2} }
\def\mat{\mathop\mathrm{Mat}\nolimits}
\def\0{\emptyset}		\def\sht{stable homotopy type}

\def\ol{\overline}	\def\dagg{\!\dagger}

\def\Si{\Sigma}	\def\Ga{\Gamma}
\def\al{\alpha}	\def\be{\beta}	\def\ga{\gamma}
\def\si{\sigma}		\def\eps{\varepsilon}

\def\rH{\mathrm H}

\def\fX{\mathbf x}	\def\bA{\mathbf A}	\def\bB{\mathbf B}
	\def\bS{\mathbf S}	

\def\sU{\mathsf U}	\def\sW{\mathsf W}
\def\sH{\mathsf H}	\def\sG{\mathsf G}

\def\mN{\mathbb N}	\def\mZ{\mathbb Z}
\def\mR{\mathbb R}	\def\mC{\mathbb C}
	\def\mP{\mathbb P}

\def\kL{\mathcal L}	\def\kI{\mathcal I}
\def\kJ{\mathcal J}	\def\kE{\mathcal E}
\def\kF{\mathcal F}	\def\kX{\mathcal X}
\def\kS{\mathcal S}	\def\kB{\mathcal B}

\title[Matrix problems and homotopy types]{Matrix problems and stable homotopy types of polyhedra}

\author{Yuriy A. Drozd}

\address{Department of Mechanics and Mathematics\\ Kyiv Taras Shevchenko University\\ 01033 Kyiv\\ Ukraine}
\email{yuriy@drozd.org}
\urladdr{drozd.org/$\sim$yuriy}

\begin{document}

\begin{abstract}
 It is a survey of the results on stable homotopy types of polyhedra of small dimensions,
 mainly obtained by H.-J.~Baues and the author \cite{bd1,bd3,bd4}. The proofs are based on the technique
 of matrix problems (bimodule categories).
 \end{abstract}

\maketitle 

\tableofcontents

\bigskip\noindent
 This paper is a survey of some recent results on stable homotopy types of polyhedra. The
 common feature of these results is that their proofs use the technique of the so called
 \emph{matrix problems}, which was mainly elaborated within framework of representation
 theory. I think that this technique is essential in homotopy theory too, and perhaps even
 in much more general setting of triangulated categories. I hope that the considerations of
 Section~3 are persuasive enough. Certainly, I could not cover
 all such results, restricting by the stable homotopy classification of polyhedra of small
 dimensions obtained in \cite{bd1,bd3,bd4,bh}. I tried to present these results in a homogeneous
 way and also replace references to rather sophisticated topological sources by simpler ones.
 The latter mainly concern some basic facts about homotopy groups of spheres, which can be
 found in \cite{hu} or \cite{tod}. I also used the book \cite{sw} as a standard source of references;
 maybe some readers will prefer \cite{sp} or \cite{co}. Most of these references are
 collected in Section~1. For the matrix problems I have chosen
 the language of \emph{bimodule categories} explained in Section~2, since it seems the simplest
 one as well as the most appropriate for applications. 

 Note that almost the same arguments that are used in Sections 5 and 6 can be applied to
 the classification of polyhedra with only $2$ non-trivial homology groups \cite{bd4}, while
 the dual arguments were applied to the spaces with only $2$ non-trivial homotopy
 groups in \cite{bd2}. Rather similar are also calculations in \cite{hen} (see also the
 Appendix by Baues and Henn to \cite{bd1}). I hope that any diligent reader of this
 survey will be able to comprehend the arguments of these papers too.

 I am extremely indebted to H.-J.~Baues, who was my co-author and my guide to the
 topological problems, and to C.~M.~Ringel, whose wonderful organising activity
 had made such a pleasant and fruitful collaboration possible. Most of our joint results
 H.-J.~Baues and I obtained during my visits to the Max-Plank-Institut f\"ur
 Mathematik, and I highly acknowledge its support.

 \section{Generalities on stable homotopy types}
 \label{s1}

 All considered spaces are supposed \emph{pathwise connected} and \emph{punctured}; we denote by $*_X$
 (or by $*$ if there can be no ambiguity) the
 marked point of the space $X$. $B^n$ and $S^{n-1}$ denote respectively the $n$-dimensional \emph{ball}
 $\setsuch{\fX\in\mR^n}{||\fX||\le1}$ and the $(n-1)$-dimensional \emph{sphere} $\setsuch{\fX\in\mR^n}{||\fX||=1}$,
 both with the marked point $(1,0,\dots,0)\,$.
 As usually, we denote by $X\vee Y$ the \emph{bouquet} (or one point union) of $X$ and $Y$, i.e. the
 factor space $X\sqcup Y$ by the relation
 $*_X=*_Y$, and identify it with $*_X\xx Y\cup X\xx*_Y\sb X\xx Y$; we denote by $X\wee Y$ the factor
 space $X\xx Y/X\vee Y$. In particular, we denote by $\Si X=S^1\wee X$ the \emph{suspension} of $X$ and by
 $\Si^nX=\underbrace{\Si\dots\Si}_{n\mbox{ \scriptsize times}}X$ its $n$-th suspension.
 The word ``\emph{polyhedron}'' is used as a synonym of ``\emph{finite \CW}.''  One can also consider bouquets of
 several spaces $\bigvee_{i=1}^sX_i$; if all of them are copies of a fixed space $X$, we denote such a bouquet by $sX$.
 
 We recall several facts on stable homotopy category of \CWs. We denote by $\hot(X,Y)$
 the set of homotopy classes of continuous maps $X\to Y$ and by $\cw$ the \emph{homotopy category} of polyhedra,
 i.e. the category whose objects are polyhedra and morphisms are homotopy classes of continuous maps.
 The suspension functor defines a natural map $\hot(X,Y)\to\hot(\Si X,\Si Y)$. Moreover, the Whitehead theorem
 \cite[Theorem 10.28]{sw} shows that the suspension functor \emph{reflects isomorphisms}: if $\Si f$ is an
 isomorphism (i.e. a homotopy equivalence), so is $f$. We set $\hos(X,Y)=\dlim_n\hot(\Si^n X,\Si^n Y)$.
 If $\al\in\hot(\Si^n X,\Si^n Y),\ \be\in\hot(\Si^m Y,\Si^m Z)$,
 one can consider the class $\Si^n\be\circ\Si^m\al\in\hot(\Si^{m+n}X,\Si^{n+m}Z)$, which is, by definition, the product
 $\be\al$ of the classes of $\al$ and $\be$ in $\hos(X,Z)$. Thus we obtain the \emph{stable homotopy category} of polyhedra
 $\cws$. Actually, if we only deal with \emph{finite \CWs}, we need not go too far, since the Freudenthal theorem
 \cite[Theorem 6.26]{sw} implies the following fact.

 \begin{proposition}\label{11}
  If $X,Y$ are of dimensions at most $d$ and $(n-1)$-connected, where $d<2n-1$, then the map $\hot(X,Y)\to\hot(\Si X,\Si Y)$ is 
 bijective.  If $d=2n-1$, this map is surjective. In particular, the map $\hot(\Si^mX,\Si^mY)\to\hos(X,Y)$
 is bijective if $m> d-2n+1$ and surjective if $m=d-2n+1$.

 Here $(n-1)$\emph{-connected} means, as usually, that $\pi_k(X)$, the $k$-th homotopy group of $X$, is trivial for
 $k\le n-1$. Thus for \emph{all} polyhedra of dimension at most $d$ the map $\hot(\Si^mX,\Si^mY)\to\hos(X,Y)$ is bijective
 if $m\ge d$ and surjective if $m=d-1$.

 Note also that the natural functor $\cw\to\cws$ reflects isomorphisms.
 \end{proposition} 

 Since we are only interested in stable homotopy classification, we identify, in what follows, polyhedra and continuous
 maps with their images in $\cws$. We denote by $\cwf$ the full subcategoy of $\cws$ consisting of all spaces $X$ with
 torsion free homology groups $\rH_i(X)=\rH_i(X,\mZ)$ for all $i$.

 Recall that any suspension $\Si^nX$ is an \emph{H-cogroup} \cite[Chapter 2]{sw}, commutative if $n\ge 2$, the category $\cws$
 is an additive category. Moreover, one can deduce from the Adams' theorem \cite[Theorem 9.21]{sw} that this category is
 actually \emph{fully additive}, i.e. every idempotent $e\in\hos(X,X)$ splits. In this case it means that there is a decomposition
 $\Si^mX\iso Y\vee Z$ for some $m$, such that $e$ comes from the map $\eps:Y\vee Z\to Y\vee Z$ with
 $\eps(y)=y$ for $y\in Y$ and $\eps(z)=*_{Y\vee Z}$ for $z\in Z$. We call a polyhedron $X$ \emph{indecomposable}
 if $X\iso Y\vee Z$ implies that either $Y$ or $Z$ are contractible (i.e. isomorphic in $\cw$ to the 1-point space).
 Actually, the category $\cws$ is a \emph{triangulated category} \cite{gm}. The suspension plays the role of shift,
 while the \emph{triangles} are the \emph{cone sequences} $X\str f Y\to Cf\to \Si X$ (and isomorphic
 ones), where $Cf=CX\cup_f Y$ is the
 \emph{cone of the map} $f$, i.e the factor space $CX\sqcup Y$ by the relation $(x,0)\sim f(x)$; $CX=X\xx I/X\xx 1$
 is the cone over the space $X$. Note that cone sqeuences coincide with \emph{cofibration sequences} in the category
 $\cws$ \cite[Proposition 8.30]{sw}. Recall that a cofibration sequence is a such one
 \begin{equation}\label{e11} 
  X \str f Y \str g Z \str h \Si X \str{\Si f} \Si Y
 \end{equation}
 that for every polyhedron $P$ the induced sequences 
 \begin{equation}\label{e12}  
  \begin{split}
 &   \hos(P,X) \str{f_*} \hos(P,Y) \str{g_*} \hos(P,Z) \str{h_*} \hos(P,\Si X) \str{\Si f_*} \hos(P,\Si Y), \\   
 &   \hos(\Si Y,P) \str{\Si f^*} \hos(\Si X,P) \str{h^*} \hos(Z,P) \str{g^*} \hos(Y,P) \str{ f^*} \hos(X,P)
 \end{split}
 \end{equation} 
 are exact. In particular, we have an exact sequence of \emph{stable homotopy groups}
 \begin{equation}\label{e13} 
    \pi^S_k(X) \str {f_*} \pi^S_k(Y) \str{g_*} \pi^S_k(Z) \str{h_*} 
 \pi^S_{k-1}(X) \str{\Si f_*} \pi^S_{k-1}(Y),
 \end{equation}
 where $\pi^S_k(X)=\dlim_m\pi_{k+m}(\Si^mX)=\hos(S^k,X)$.
 Certainly, one can prolong the sequences \eqref{e12} and \eqref{e13} into infinite exact
 sequences just taking further suspensions.
 
 Every \CW\ is obtained by \emph{attaching cells}. Namely, if $X^n$ is the $n$-th skeleton of $X$, then there is
 a bouquet of balls $B=mB^{n+1}$ and a map $f:mS^n\to X^n$ such that $X^{n+1}$ is isomorphic to the cone
 of $f$, i.e. to the space $X^n\cup_fB$.  It gives cofibration sequences like \eqref{e11} and exact sequences
 like \eqref{e12} and \eqref{e13}.

 We denote by $\cw^k_n$ the full subcategory of $\cw$ formed by $(n-1)$-connected $(n+k)$-dimensional polyhedra and by
 $\cwf^k_n$ the full subcategory of $\cw^k_n$ formed by the polyhedra $X$ with torsion free homology groups $\rH_i(X)$
 for all $i$. Proposition \ref{11} together with the fact that every map of \CWs\ is homotopic to
 a cell map, also implies the following result.

  \begin{proposition}\label{12}
   The suspension functor $\Si$ induces equivalences $\cw^k_n\to\cw^k_{n+1}$ for all $n>k+1$. Moreover, if $n=k+1$,
 the suspension functor $\Si:\cw^k_n\to\cw^k_{n+1}$ is a \emph{full representation equivalence}, i.e. it is full, dense
 and reflects isomorphisms. (\emph{Dense} means that every object
 from $\cw^k_{n+1}$ is isomorphic (i.e. homotopy equivalent) to $\Si X$ for some $X\in\cw^k_n$.)
 \end{proposition} 

 Therefore, setting $\cw^k=\cw^k_{k+2}\iso\cw^k_n$ for $n>k+1$, we can consider it as a full subcategory of
 $\cws$. The same is valid for $\cwf^k=\cwf^k_{k+2}$. Note also that $\cw^k_n$ naturally embeds into $\cw^k_{n+1}$.
 It leads to the following notion \cite{ba1}.

 \begin{definition}\label{13} 
 An \emph{atom} is an indecomposable polyhedron $X\in\cw^k_{k+1}$ not belonging to the image of $\cw^k_k$.
 A \emph{suspended atom} is a polyhedron $\Si^mX$, where $X$ is an atom.
 \end{definition} 

 Then we have an obvious corollary.

  \begin{corollary}\label{14}
  Every object from $\cw^k_n$ with $n\ge k+1$ is isomorphic (i.e. homotopy equivalent) to a bouquet $\bigvee_{i=1}^sX_i$,
 where $X_i$ are suspended atoms. Moreover, any suspended atom is indecomposable (thus indecomposable objects are just
 suspended atoms). 
 \end{corollary} 

 Note that the decomposition in Corollary \ref{14} is, in general, not unique \cite{fr}. That is why an important question is the structure of
 the \emph{Grothendieck group} $K_0(\cw^k)$. By definition, it is the group generated by the isomorphism classes $[X]$ of
 polyhedra from $\cw^k$ subject to the relations $[X\vee Y]=[X]+[Y]$ for all possible $X,Y$. The following results of Freyd \cite{fr,co}
 describe the structure of this group.

  \begin{definition}\label{15}
  \begin{enumerate}
\item  Two polyhedra $X,Y\in\cw^k$ are said to be \emph{congruent} if there is a polyhedron $Z\in\cw^k$ such that
 $X\vee Z\iso Y\vee Z$ (in $\cw^k$).
 \item	A polyhedron $X\in\cw^k$ is said to be $p$\emph{-primary} for some prime number $p$ if there is a bouquet of spheres
 $B$ such that the map $p^m1_X:X\to X$ can be factored through $B$, i.e. there is a commutative diagram 
 $$
   \xymatrix{	X \ar[rr]^{p^m1_X}\ar[dr] && X \\
		& B \ar[ur] } 
 $$
\end{enumerate}	
 \end{definition} 

  \begin{theorem}[Freyd]\label{16}
  The group $K_0(\cw^k)$ (respectively $K_0(\cwf^k)\,$) is a free abelian group with a basis formed by the
 congruence classes of  $p$-primary suspended atoms from $\cw^k$ (respectively from $\cwf^k$) for all
 prime numbers $p\in\mN$.
 \end{theorem} 

 Therefore, if we know the ``place'' of every atom class $[X]$ in $K_0(\cw^k)$ or $K_0(\cwf^k)$, i.e.
 its presentation as a linear combination of classes of $p$-primary suspended atoms, we can deduce herefrom
 all decomosition rules for $\cw^k$ or $\cwf^k$.

 \section{Bimodule categories}
 \label{bc}
 
 We also recall main notions concerning \emph{bimodule categories} \cite{d1,d2}. Let $\bA,\bB$ be two fully additive
 categories. An \emph{$\bA$-$\bB$-bimodule} is, by definition, a biadditive bifunctor $\sU:\bA^\circ\xx\bB\to\Ab$.
 As usually, given an element $u\in\sU(A,B)$ and morphisms $\al\in\bA(A',A),\ \be\in\bB(B,B')$, we write $\be u\al$
 instead of $\sU(\al,\be)u$.
 Given such a functor, we define the \emph{bimodule category $\El(\sU)$} (or the category of \emph{elements of the bimodule}
 $\sU$, or the category of \emph{matrices over} $\sU$) as follows.
 \begin{itemize}
\item 	The set of \emph{objects} of $\El(\sU)$ is the disjoint union 
 $$\ob\El(\sU)=\bigsqcup_{\substack{A\in\ob\bA\\B\in\ob\bB}}\sU(A,B).$$
 \item	A \emph{morphism} from $u\in\sU(A,B)$ to $u'\in\sU(A',B')$ is a pair $(\al,\be)$ of morphisms $\al\in\bA(A,A'),\,
 \be\in\bB(B,B')$ such that $u'\al=\be u$ in $\sU(A,B')$. 
 \item	The product $(\al',\be')(\al,\be)$ is defined as the pair $(\al'\al,\be'\be)$.
 \end{itemize}
 Obviously, $\El(\sU)$ is again a fully additive category. 

 Suppose that $\ob\bA\supset\set{\lst An},\ \ob\bB\supset\set{\lst Bm}$ such that
 every object $A\in\ob\bA$ \,($B\in\ob\bB$) decomposes as $A\iso\b+_{i=1}^nk_iA_i$ (respectively, $B\iso\b+_{i=1}^ml_iB_i$).
 Then $\bA^{\circ}$ (respectively, $\bB$) is equivalent to the category of finitely generated projective
 right (left) modules over the ring of matrices
 $(a_{ij})_{n\xx n}$ with $a_{ij}\in \bA(A_j,A_i)$ (respectively, $(b_{ij})_{m\xx m}$ with $b_{ij}\in\bB(B_j,B_i)$). We denote these
 rings respectively by $|\bA|$ and $|\bB|$. We also denote by $|\sU|$ the $|\bA|\mbox{-}|\bB|$-bimodule consisting of matrices
 $(u_{ij})_{m\xx n}$, where $u_{ij}\in\sU(A_j,B_i)$. Then $\sU(A,B)$, where $A,B$ are, respectively, a projective right
 $|\bA|$-module and a projective left $|\bB|$-module, can be identified with $A\*_{|\bA|}|\sU|\*_{|\bB|}B$. Elements from this set
 are usually considered as block matrices $(U_{ij})_{m\xx n}$, where the block $U_{ij}$ is of size $l_i\xx k_j$ with entries from
 $\sU(A_j,B_i)$. To form a direct sum of such elements, one has to write direct sums of the corresponding blocks at each place.
 Certainly, some of these blocks can be ``empty,''  if $k_j=0$ or $l_i=0$. An empty block is indecomposable \iff it is of size $0\xx 1$
 (in $\sU(A_j,0)$) or $1\xx 0$ (in $\sU(0,B_i)\,$); we denote it respectively by $\0^j$ or by $\0_i$.

 In many cases the rings $|\bA|$ and $|\bB|$ can be identified with \emph{tiled subrings} of rings of integer matrices. Here a 
 \emph{tiled subring} in $\mat(n,\mZ)$ is given by an integer matrix $(d_{ij})_{n\xx n}$ such that $d_{ii}=1$ and $d_{ik}|d_{ij}d_{jk}$
 for all $i,j,k$; the corresponding ring consists of all matrices $(a_{ij})$ such that $d_{ij}|a_{ij}$ for all $i,j$
 (especially $a_{ij}=0$ if $d_{ij}=0$).
 
  \begin{example}\label{bc-1}
  Let $\bA\subset\mat(2,\mZ)$ be the tiled ring given by the matrix 
 $$ 
    \begin{pmatrix}
  1 & 12 \\ 0 & 1
 \end{pmatrix} ,
 $$ 
   $\sU$ be the set of $2\xx2$-matrices $(u_{ij})$ with $u_{ij}\in \mZ/24$ if $i=1,j=2$, $u_{ij}\in\mZ/2$ otherwise. We define $\sU$
 as an $\bA\mbox{-}\bA$-bimodule setting 
 \begin{align*}
    \begin{pmatrix} a&12b\\ 0&c\end{pmatrix}  \begin{pmatrix} u_1&u_2\\u_3&u_4\end{pmatrix}=&  \begin{pmatrix} au_1+bu_3\ &
 au_2+12bu_4\\cu_3&cu_4 \end{pmatrix} ;\\
 \begin{pmatrix} u_1&u_2\\u_3&u_4 \end{pmatrix}  \begin{pmatrix} a&12b\\ 0&c\end{pmatrix} =&
 \begin{pmatrix} au_1 &\ cu_2+12bu_1\\ au_3 &cu_4+bu_3\end{pmatrix} . 
 \end{align*} 
 If we need to indicate this action, we write 
 $$ 
       \begin{pmatrix}
  1 & 12^* \\ 0 & 1
 \end{pmatrix}\ \text{ and }\  \begin{pmatrix}
  \mZ/2 &\ \mZ/24 \ \\\ \mZ/2^* & \mZ/2 
 \end{pmatrix} 
 $$ 
 for the matrix defining the ring $\bA$ and for the bimodule $\sU$. Thus the multiplications of the elements marked
 by stars is given by the \emph{$\ast$-rule}:
 \begin{equation}\label{e21} 
  (12a^*)\cdot (u\mod2^*)=ac\mod2.
 \end{equation} 
 \end{example} 

  \begin{example}\label{bc-1a} 
 In the classification of torsion free atoms below the following bimodule plays the crucial role.
 We consider the tiled rings $\bA_2\subset\mat(2,\mZ)$
 and $\bB_2\subset\mat(7,\mZ)$ given respectively by the matrices  
 $$ 
  \begin{pmatrix}
  1 & 2 & 2 & 12 & 24 & 12 & 24 \\ 1 & 1 & 1 & 12 & 24 & 6 & 24 \\ 1 & 2 & 1 & 12 & 24 & 12 & 24 \\
 0&0&0& 1&2&12^*&12\\ 0&0&0& 1&1&12&6\\ 0&0&0& 0&0&1&1& \\ 0&0&0& 0&0&0&1
 \end{pmatrix}\ \text{ and }\    \begin{pmatrix}
   1 & 12^* \\ 0 & 1
 \end{pmatrix}  .
 $$ 
 The $\bA_2\mbox{-}\bB_2$-bimodule $\sU_2$ is defined as the set of matrices of the form 
 $$ 
    \begin{pmatrix}
 \ \mZ/24 & 0\\ \mZ/12 &0 \\ \mZ/12 &0\\ \mZ/2 &\ \mZ/24 \ \\ 0 &\mZ/12 \\ \mZ/2^* &\mZ/2 \\  0 &\mZ/2
 \end{pmatrix} .
 $$ 
 The multiplication in $\sU_2$ is given by the natural matrix multiplication, but taking into account the $\ast$-rule
 \eqref{e21}.   
 
 We shall use the following description of indecomposable elements in $\El(\sU_2)$.
 Set $I_1=\set{1,2,3,4,6},\,I_2=\set{4,5,6,7}$, $V=\setsuch{v\in\mN}{1\le v\le6}$,
 $V_1=\setsuch{v\in\mN}{1\le v\le 12}$, $V_2=\set{1,2,3}$.
 \end{example}

 \begin{theorem}\label{bc-2}
  A complete list $\kL_2$ of non-isomorphic indecomposable objects from $\El(\sU_2)$ consists of
\begin{itemize}
 \item 	empty objects $\0^j\ (j=1,2)$ and $\0_i\ (1\le i\le 7)$;
 \item	objects $v^j_i\in\sU(A_j,B_i)\ (j=1,2;\,i\in I_j;\ v\in V_1$ if $i=1;\ v=1$ if $i=6,7$ or $(ij)=(14);\ v\in V$ otherwise); 
 \item	objects $v^j_{il}=\ds\binom {v^j_i}{1^j_l}\ (j=1,2;\ i=1,2,3,\,l=4,6$ if $j=1;\ i=4,5,\,l=6,7$ if $j=2$; if $(il)=(26)$ or $(57)$
 then $v\in V_2$; otherwise $v\in V$);
 \item	objects $v_{44}=(1^1_4 \ v^2_4)$ with $v\in V$;
 \item	objects $v_{4l}=\ds \begin{pmatrix}
  1^1_4 & v^2_4 \\ 0& 1^2_l
 \end{pmatrix} $ with $l=6,7$ and $v\in V$;
 \item	objects $v_iw_{44}=\ds \begin{pmatrix}
   v^1_i &0 \\ 1^1_4 & w^2_4 
 \end{pmatrix} $ with $i=1,2,3$ and $v,w\in V$;
 \item	objects $v_iw_{4l}=\ds \begin{pmatrix}
  v^1_i & 0\\ 1^1_4 & w^2_4 \\ 0& 1^2_l
 \end{pmatrix} $ with $i=1,2,3,\,l=6,7$ and $v,w\in V$.
\end{itemize}

\medskip\noindent
 Here the indices define the block containing the corresponding element.  
 \end{theorem} 
 \begin{proof}
   Decompose $\sU$ into $2$-primary and $3$-primary parts. Since for every two matrices $M_2,M_3\in\GL(n,\mZ)$ there is a
 matrix $M\in\GL(n,\mZ)$ such that $M\equiv M_2\mod2$ and $M\equiv M_3\mod3$, we can consider $2$-primary part
 and $3$-primary part separately. Note that in the $3$-primary part the blocks $u^1_4,u^1_6,u^2_6$ and $u^2_7$ vanish, while
 the other non-zero blocks of $u\in\ob(\sU_2)$ are with entries from $\mZ/3$ and there are no restrictions on elementary transformation
 of the matrix $u$. Thus every element in the $3$-primary part is a direct sum of elements $1^j_i$ with $j=1,i=1,2,3$ or $j=2,i=4,5$.

 For elements $u,u'$ of the $2$-primary part write $u<u'$ if $u'=ua$ for some non-invertible $a\in\bA_2$. Then we have the following
 inequalities:
 \begin{align*}
  & 1^1_1<1^1_3<1^1_2<2^1_1<2^1_3<2^1_2<4^1_1,\\
  & 1^1_6<1^1_4<4^1_1 \ \text{ and }\ 1^1_6<2^1_2; \\
 &  1^2_4<1^2_5<2^2_4<2^2_5<4^2_4,\\
 &  1^2_7<1^2_6<4^2_4 \ \text{ and }\ 1^2_7<2^2_5.
 \end{align*} 
  Using them, one can easily decompose the parts 
 $$ 
  \tilde u^1=  \begin{pmatrix} u^1_1\\ u^1_2\\ u^1_3\end{pmatrix} \ \text{ and }\ 
 \tilde u^2= \begin{pmatrix}u^2_4\\ u^2_5 \end{pmatrix}
 $$ 
 into a direct sum of empty and $1\xx 1$ matrices. Now we obtain a column splitting of the remaining matrices, and with respect
 to the transformation that do not change $\tilde u^1$ and $\tilde u^2$, these columns are linearly ordered. Therefore, we can also
 split them into empty and $1\xx 1$ blocks. Together with $\tilde u^1$ and $\tilde u^2$, it splits the whole matrix $u$ into a direct sum
 of matrices of the forms from the list $\kL_2$, where $v,w$ are powers of $2$. Adding $3$-primary parts, we get the result.
 \end{proof} 
 
 \begin{example}\label{bc-1b} 
 Consider the idempotents $e=\sum_{i\in I_1}e_{ii}\in\bA_2$ and $e'=e_{11}\in\bB$. Set $\bA_1=e\bA e,\,\bB_1=e'\bB_2e'
\iso\mZ$ and $\sU_1=e'\sU_2 e$. Then $\sU_1$ is an $\bA_1\mbox{-}\bB_1$-bimodule; elements from $\El(\sU_1)$
 can be identified with those from $\El(\sU_2)$ having no second column and fifth row. Hence we get the following
 result.  
 \end{example}

 \begin{corollary}\label{bc-3}
    A complete list $\kL_1$ of non-isomorphic indecomposable objects from $\El(\sU_1)$ consists of
\begin{itemize}
 \item 	empty objects $\0_i\ (i\in I_1)$;
 \item	objects $v_i\ (i\in I_1;\ v\in V_1$ if $i=1,\ v\in V$ if $i=2,3,\ v=1$ if $i=4,6$); 
 \item	objects $v_{il}=\ds\binom {v_i}{1_l}\ (i=1,2,3,\,l=4,6$; if $(il)=(26)$ then $v\in V_2$, otherwise $v\in V$).
 \end{itemize}

\medskip\noindent
 Here the indices show the blocks where the corresponding elements are placed.
 \end{corollary} 

 \section{Bimodules and homotopy types}
 \label{bh}

  Bimodule categories arise in the following situation. Let $\bA$ and $\bB$ are two fully additive subcategories of
 the category $\Hos$. We denote by $\bA\dagg\bB$ the full subcategory of $\Hos$ consisting of all objects $X$
 isomorphic (in $\Hos$) to the cones of morphisms $f:A\to B$ with $A\in\bA,\,B\in\bB$, or, the same,
 such that there is a cofibration sequence
 \begin{equation}\label{bh-e1} 
   A \str f B \str{g} X \str{h} \Si A, 
 \end{equation}
 where $A\in\bA,\,B\in\bB$. Consider the $\bA$-$\bB$-bimodule $\sH$, which is the restriction on $\bA^\circ\xx\bB$
 of the ``regular'' $\Hos$-$\Hos$-bimodule $\hos$. If $f\in\hos(A,B)$ is an element of $\sH$, it gives rise to an exact
 sequence like \eqref{bh-e1} with $X=Cf$. Moreover, since this sequence is a cofibration one, every morphism
 $(\al,\be):f\to f'$, where $f'\in\Hos(A',B')$ induces a morphism $\ga:X\to X'$, where $X'=Cf'$, such that the diagram
 \begin{equation}\label{bh-e2} 
   \begin{CD}
  A @>f>> B @>g>> X @>h>> \Si A @>\Si f>> \Si B\\
 @V\al VV  @V\be VV  @VV\ga V  @VV\Si\al V @VV\Si\be V\\
  A' @>>f'> B'@>>g'> X' @>>h'> \Si A' @>>\Si f'> \Si B'
 \end{CD}
 \end{equation}
 commutes. In what follows we suppose that the categories $\bA$ and $\bB$ satisfy the following condition:
 \begin{equation}\label{bh-e3} 
  \hos(B,\Si A)=0 \quad\text{for all}\quad A\in\bA,\ B\in\bB. 
  \end{equation}
 In this situation, given a morphism $\ga:X\to X'$, we have that $h'\ga g=0$, hence $\ga g=g'\be$ for some
 $\be:B\to B'$. Moreover, since the sequence 
 $$ 
  B \str g X \str h \Si A \str{\Si f} \Si B  
 $$ 
 is cofibration as well, and $\Si:\Hos(A,B)\to\Hos(\Si A,\Si B)$ is a bijection, there is a morphism
 $\al:A\to A'$, which makes the diagram \eqref{bh-e2} commutative.

 Note that neither $\ga$ is uniquely determined by $(\al,\be)$, nor $(\al,\be)$ is uniquely restored from $\ga$. 
 Nevertheless, we can control this non-uniqueness. Namely, if both $\ga$ and $\ga'$ fit the diagram \eqref{bh-e2}
 for given $(\al,\be)$, their difference $\ol\ga=\ga-\ga'$ fits an analogous diagram with $\al=\be=0$. The equality
 $\ol\ga g=0$ implies that $\ol\ga=\si h$ for some $\si:\Si A\to X'$, and the equality $h'\ol\ga=0$ implies that
 $\ol\ga=g'\tau$ for some $\tau:X\to B$. On the contrary, if $\ol\ga=\si h=g'\tau$ for \emph{some}
 morphisms $\si:X\to \Si Y\,\tau:X\to Z$, where $Y\in\bA,\,Z\in\bB$, the condition \eqref{bh-e3} implies that
 $\ol\ga g=h'\ol\ga=0$, so $\ol\ga$ fits the diagram \eqref{bh-e2} with $\al=\be=0$. 

 Fix now $\ga$, and let both $(\al,\be)$ and $(\al',\be')$ fit \eqref{bh-e2} for this choice of $\ga$. Then the pair
 $(\ol\al,\ol\be)$, where $\ol\al=\al-\al',\,\ol\be=\be-\be'$, fits \eqref{bh-e2} for $\ga=0$. The equality $g'\ol\be=0$
 implies that $\ol\be=f'\si$ for some $\si:B\to A'$, and the equality $(\Si\al)h=0$ implies that $\Si\al=\Si\tau\Si f$,
 or $\al=\tau f$ for some $\tau:B\to S$. On the contrary, if $(\ol\al,\ol\be):f\to f'$ is such that $\ol\be=f'\si$
 and $\ol\al=\tau f$ with $\si,\tau:B\to A'$, then $g'\be=(\Si\al) h=0$, hence this pair fits \eqref{bh-e2} with $\ga=0$. 

 Summarizing these considerations, we get the following statement.

 \begin{theorem}\label{bh-1} 
  Let $\bA,\bB$ be fully additive subcategories of $\Hos$ satisfying the condition \eqref{bh-e3}, $\bA\dagg\bB$ be
 the full subcategory of $\Hos$ consisting of all spaces such that there is a cofibration \eqref{bh-e1} with $A\in\bA$,
 $B\in\bB$. Denote by $\sH$ the bimodule $\hos$ considered as $\bA$-$\bB$-bimodule, by $\kI$ the ideal in
 $\bA\dagg\bB$ consisting of all morphisms $\ga:X\to X'$ that factor both through an object from $\Si\bA$
 and through an object from $\bB$, and by $\kJ$ the ideal in $\El(\sH)$ consisting of all morphisms
 $(\al,\be):f\to f'$ such that $\be$ factors through $f'$ and $\al$ factors through $f$.
 Then the factor categories $\El(\sH)/\kJ$ and $\bA\dagg\bB/\kI$ are equivalent; an equivalence is induced
 by the maps $f\mapsto Cf$ and $(\al,\be)\mapsto\ga$, where $\ga$ fits a commutative diagram \eqref{bh-e2}.
 Moreover, $\kI^2=0$, thus the functor $\bA\dagg\bB\to\bA\dagg\bB/\kI$ reflects isomorphisms.
 \end{theorem} 
 \begin{proof} 
 We only have to check the last statement. But if $\ga:X\to X'$ factors as $X\str\tau B\str{g'}X'$ and
 $\ga':X'\to X''$ factors as $X'\str h' \Si A\str\si X''$, where $A\in\bA,\,B\in\bB$, then $\ga'\ga=0$,
 since $h'g:B\to\Si A$ and $\hos(B,\Si A)=0$.
 \end{proof} 

   \begin{corollary}\label{bh-2} 
   In the situation of Theorem~\ref{bh-1}, suppose that $\hos(B,A)=0$ for each $A\in\bA,\,B\in\bB$. Then
 $\El(\sH)\iso\bA\dagg\bB/\kI$. Moreover, the functor $\bA\dagg\bB\to\El(\sH)$ is a \emph{representation equivalence},
 i.e. it is dense, preserves indecomposables and reflects isomorphisms.
 \end{corollary}

 Note also that any \emph{isomorpism} $f:A\str\sim B$ is a zero object in $\El(\sH)/\kJ$, since its identity
 map $(1_A,1_B)$ can be presented as $(f^{-1}f,ff^{-1})$. Obviously, the corresponding object from
 $\bA\dagg\bB$ is zero (i.e. contractible) too.

 \section{Small dimensions}
 \label{sd} 

 We now use Theorem~\ref{bh-1} to describe \sht s of atoms of dimensions at most 5, or, the same,
 indecomposable objects in the categories $\cw^1_2$ and $\cw^2_3$.

  \begin{example}\label{sd-1}
  It is well known that $\pi_n(S^n)=\mZ$ (freely generated by the identity map). It allows easily to describe 
  atoms in $\cw^1_2$. Such an atom $X$ is (stably!) of the form $Cf$ for some map $f:mS^2\to nS^2$.
 Since $\hos(S^n,S^{n+1})=0$, Theorem~\ref{bh-1} can be applied. The map $f$ is given by an integer
 matrix. Using automorphisms of $mS^2$ and $nS^2$, we can transform it to a diagonal form.
 Hence, indecomposable gluings can only be if $m=n=1$; thus $f=q1_{S^2}$. One can see that such a gluing
 is indecomposable \iff $q$ is a power of a prime number. The corresponding atom $S^2\cup_qB^3$ will be
 denoted by $M(q)$ and called \emph{Moore atom}. It occurs in a cofibration sequence 
 $$ 
  S^2 \str q S^2 \xarr{g(q)} M(q) \xarr{h(q)} S^3 \str q S^3.
 $$ 

 For the next section we need more information about $2$-primary Moore atoms. We denote $M_t=M(2^t)$ and write
 $g_t,h_t$ instead of $g(2^t),h(2^t)$. These atoms can be included into 
 the following commutative ``octahedral'' diagram \cite{gm}, where $t=r+s$:
 \begin{equation}\label{sd-e1} 
  \xy<0em,2.5em>
  \xymatrix@R=1em{
	& S^2 \ar[rr]^{g_r} \ar[dd]^{\,2^s} && M_r \ar[dr]^{h_r} \ar[dd]_{k_{tr}} \\
	S^2 \ar[ur]^{2^r} \ar[dr]_{2^t} &&&& S^3	\\
	& S^2 \ar[rr]_{g_t} \ar[dr]_{g_s} && M_t \ar[ur]_{h_t}\ar[dl]^{k_{st}} \\
	&& M_s
  } \endxy
 \end{equation} 
 Moreover, in this diagram $h_sk_{st}=2^rh_t$.

 The exact sequence \eqref{e13} is here of the form 
 $$ 
   \pi^S_k(S^2)\str q\pi^S_k(S^2)\larr\pi^S_k(M(q))\larr \pi^S_k(S^3)\str q \pi^S_k(S^3), 
 $$ 
 which gives the values of stable homotopy groups of  the spaces $M(q)$ shown in Table~\ref{tab1} below.
 \begin{table}[!ht]
 \caption{}\label{tab1}
$$  \begin{array}{|c|c|c|c|}
  \hline 
  k &	  2 & 3 & 4  \\
 \hline
 \ \pi^S_k(M(q)),\ q \text{ odd}\  &\ \mZ/q\ & 0 & 0  	\\
 \hline
 \ \pi^S_k(M_t),\ t>1 \  & \mZ/q &\ \mZ/2\ &\ \mZ/2\+\mZ/2\ 	\\
\hline
  \pi^S_k(M_1) & \mZ/2 & \mZ/2 & \mZ/4 \\
\hline
 \end{array} 
 $$
\end{table}
 (By the way, this table implies that all Moore atoms are pairwise non-isomorphic.)
 Actually, the only non-trivial is the group $\pi^S_4(M_1)$. It can be obtained as $\pi_6(\Si^2M_1)$,
 which is isomorphic to the $2$-primary component of $\pi_6(S^3)=\mZ/12$ (cf. \cite[Lemma XI.10.2]{hu}).
 To prove that the sequence
 $$ 
  0\larr \pi^S_4(S^2)=\mZ/2 \larr \pi^S_4(M_t) \larr \pi^S_4(S^3)=\mZ/2 \larr 0  
 $$ 
 splits if $t>1$, it is enough to consider the commutative diagram 
 \begin{equation}\label{sd-e2} 
 \xy<0pt,1.5em>
  \xymatrix{
   0 \ar[r] & \pi^S_4(S^2) \ar[r] \ar[d]_{0\,} & \pi^S_4(M_1) \ar[r]\ar[d] & \pi^S_4(S_3) \ar@{=}[d] \ar[r]& 0 \\
    0 \ar[r] & \pi^S_4(S^2) \ar[r]  & \pi^S_4(M_t) \ar[r] & \pi^S_4(S_3) \ar[r] &0 
 }\endxy
 \end{equation}
 arising from the diagram \eqref{sd-e1} with $r=1$. It shows that
 the second row of this diagram is the pushdown of the first one along the zero map; thus it splits.
 \end{example} 

  \begin{example}\label{sd-2}
  Now we are able to describe atoms in $\cw^2_3$. They are cones $Cf$ for some $f:mS^4\to Y$
 with $2$-connected $Y$ of dimension $4$. Again $\Hos(Y,S^5)=0$, so Theorem~\ref{bh-1} can be applied. 
 Example~\ref{sd-1} shows that $Y$ is a bouquet of spheres $S^3,S^4$ and suspended Moore atoms
 $\Si M(q)$. Note that $\pi^S_4(Y)=\pi_4(Y)$ for every $Y$; in particular $\pi_4(S^4)=\mZ,\ \pi_4(S^3)=\mZ/2$
 (generated by the suspended Hopf map $\eta_1=\Si\eta$; $\eta:S^3\to S^2\iso\mC\mP^1$ is given by the
 rule $\eta(a,b)=(a:b)$, where $(a,b)\in\mC^2$ are such that $|a|^2+|b|^2=1$) and 
$$
\pi_4(\Si M(q))=\pi^S_3(M(q))=\begin{cases} \mZ/2 &\mbox{if }\ q=2^r \\
	0 &\mbox{otherwise}. \end{cases}
 $$
 The Hopf map $\eta_2=\Si^2\eta:S^4\to S^3$ and the inclusion $j:S^2\to M(q)$ give rise to an epimorphism
 $\eta_*:\pi_4(S_4)\to \pi_4(S^3)$ and to an isomorphism $j_*:\pi_4(S^3)\to\pi_4(\Si M_r)$, where $M_r=M(2^r)$.
 Moreover, if $t>r$, there is a map $M(2^r)\to M(2^t)$ that induces an isomorphism $\pi_4(M_r)\to\pi_4(M_t)$. 
 If $Y=s_4S^4\vee s_3S^3\vee(\bigvee_{r=1}^\8 m_rM_r)$, a map $f:mS^4\to Y$ can be given by a matrix of the form 
 $$ 
    \begin{pmatrix}
  F_4 & F_3 & G_1 & G_2 & \dots
 \end{pmatrix} ,
 $$ 
 where $F_i$ is of size $m\xx s_i$ with entries from $\pi_4(S^i)$; $G_r$ is of size $m\xx m_r$ with entries from
 $\pi_4(\Si M_r)$ (some of these matrices can be ``empty,'' containing no columns).
 Using automorphisms of $Y$ and $B$, one can easily transform this matrix to the shape where
 there is at most two non-zero elements in every row (if two, one of them necessarily in the matrix $F_4$ and even) and at most 
 one non-zero element in every column, as shown below:  
 $$ 
    \begin{array}{|c|c|c|}
 \hline
  \ F_4\ &\ F_3\ &\ G_r\ \\
 \hline
  q &  & \\
 \hline
   & \eta & \\
\hline
 2^t &\eta &\\
 \hline
 &&\eta \\
\hline
 2^t & &\eta \\
\hline
 \end{array} 
 $$ 
 Thus $X$ decomposes into a bouquet of the spaces $\Si^2M(q)$ (which are not atoms,
 but suspended atoms), spheres and the spaces $C(\eta),\,C(\eta 2^t)$ and $C(2^r\eta 2^t)$,
 which are gluings of the following forms:  

\medskip
 $$ 
   
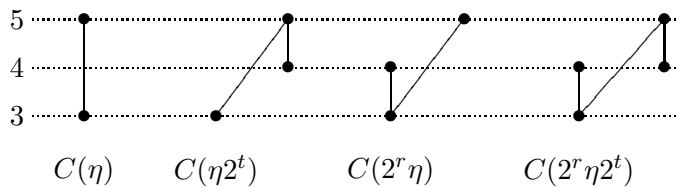
\xymatrix@=1ex{
	 {5} \ar@{.} [rrrrrrrrrrr] & \cel \ar@{-} [dd] &&& \cel \ar@{-}[ddl] &&& 
	\cel \ar@{-}[ddl]  &&& \cel \ar@{-}[ddl] \ar@{-}[d] & {} \\
        {4} \ar@{.} [rrrrrrrrrrr] &&&& \cel \ar@{-}[u] && \cel \ar@{-}[d]
	 &&& \cel\ar@{-}[d] &\cel&{} \\
        {3} \ar@{.} [rrrrrrrrrrr] & \cel && \cel &&& \cel &&& \cel&&  {} 	\\
	& {C(\eta)}  && {C(\eta 2^t)} &&& {C(2^r\eta)} &&& {C(2^r\eta 2^t)} 	} 
 $$ 

\medskip\noindent
 Here, following Baues, we denote the cells by bullets and the attaching maps by lines; the word in brackets shows which
 maps are chosen to attach bigger cells to smaller ones. We do not show the fixed point, which coincide here with $X^2$
 (since $X$ is $2$-connected); thus the lowest bullets actually describe spheres, not balls. These polyhedra are called 
 \emph{Chang atoms}. Again one can check that all of them are pairwise non-isomorphic.
 \end{example}

 Thus we have proved the following classical result.

 \begin{theorem}[Whitehead, \cite{w1}]\label{bh-5}
  The atoms of dimension at most $5$ are:
 \begin{itemize}
\item 	sphere $S^1$ (of dimension $1$);
 \item	Moore atoms $M(q)$, where $q=p^r$, $p$ is a prime number (of dimension $3$);
 \item	Chang atoms $C(\eta),\,C(\eta 2^r),\,C(2^r\eta 2^t)$ (of dimension $5$).
\end{itemize}
 \end{theorem} 

 In what follows, we often use suspended Moore and Chang atoms. We shall denote them by the same symbols
 but indicating the dimension. Thus $M^d(q)=\Si^{d-3}M(q)$ and $C^d(w)=\Si^{d-5}C(w)$ for $w\in\set{\eta,2^r\eta,
 \eta2^r,2^r\eta2^t}$; in particular, $M(q)=M^3(q)$ and $C(w)=C^5(w)$.
 The same agreement will also be used for other atoms constructed below.

 \section{Dimension 7}
 \label{d7} 

 We shall now consider the category $\cw^3$. Its object actually come from $\cw^3_4$, so we have to classify atoms
 of dimension $7$. Such an atom $X$ is $3$-connected, so we may suppose that $X^3=*$. Set $B=X^5$, then $X/B$
 only has cells of dimensions $6$ and $7$. Therefore $X\in \Si^3\cw^1\dagg\Si^2\cw^1\iso\Si\cw^1\dagg\cw^1$. Consider
 the bifunctor $\sW(A,B)=\hos(\Si A,B)$ restricted to the category $\cw^1$. Since, obviously, $\Hos(B,\Si^2A)=0$
 for $A,B\in\cw^1$, we can apply Theorem~\ref{bh-1}. So we first classify indecomposable elements of the
 bimodule category $\El(\sW)$.
 
 Indecomposable objects of the category $\cw^1$ are spheres $S^2,S^3$ and Moore atoms $M(q)\ (q=p^r,\ r$ prime). 
 If $q$ is odd, one easily sees that $\sW(A,M(q))=0$ for all $A$, so we may only consider the spaces $M_r=M(2^r)$. 

 From the cofibration sequence
 $$
 S^2\str{g_r} S^2\to M_r\str{h_r} S^3\to S^3
 $$
 and the diagram \eqref{bh-e1}, we get the values of the $\hos$-groups shown in Table~\ref{tab2}.

 \begin{table}[!ht]
 \caption{}  \label{tab2} 
 $$ 
    \begin{array}{|c|c|c|c|c|}
 \hline
  & S^2 & S^3 & M_1 &\ M_r\ (r>1)\ \\
 \hline
  S^2 &  \mZ &\ \mZ/2\ & \ \mZ/2\  & \mZ/2\\
 \hline
 S^3 & 0 & \mZ & \mZ/2 &\ \mZ/2^r\ \\
 \hline
 M_1 &\ \mZ/2\ &\ \mZ/2\ & \mZ/4&  T_{1r}\\
 \hline
\ M_t\ (t>1)\ &\ \mZ/2^t\  &\mZ/2 & T_{1t} & T_{tr} \\
 \hline
 \end{array} 
 $$ 
 \end{table}
\noindent
 Here $T_{tr}$ denotes the set of matrices  $\ds\binom{a\ b}{0\ a}$ with $a\in2^m\mZ/2^m,\,b\in\mZ/2$,
 where $m=\min(r,t)$. The equality $\hos(M_1,M_1)=\mZ/4$ follows from the fact that this ring acts on
 $\pi^S_4(M_1)=\mZ/4$, so  $2\hos(M_1,M_1)\ne0$. The diagram \eqref{sd-e1} implies that the sequence
  \begin{multline*}  
  0\larr \hos(S^3,M_t) \larr \hos(M_r,M_t) \larr\\
 \larr \ker\{\hos(S^2,M_t)\str{2^r}\hos(S^2,M_t)\}\larr 0  
 \end{multline*} 
 splits if $\min(r,t)>1$. The generator of the subgroup of diagonal matrices in $T_{tr}$ is $k_{tr}$,
 while the matrix $\ds\binom{0\ 1}{0\ 0}$ corresponds to the morphism $g_t\eta h_r$.

 Analogous calculations, using Table~\ref{tab1} of the preceding section and the diagram~\ref{sd-e1},
 produce the following Table~\ref{tab3} for the values of the functor $\hos(\Si A,B)$.
 \begin{table}[!ht]
  \caption{}\label{tab3}
 $$ 
   \begin{array}{|c|c|c|c|c|}
 \hline
 & S^2 & S^3 & M_1 &\ M_r\ (r>1)\ \\
\hline
 S^2 & \mZ/2 & \mZ/2 & \mZ/4 &\ \mZ/2\+\mZ/2\ \\
\hline
 S^3 & \mZ &\ \mZ/2\ & \mZ/2 & \mZ/2 \\ 
  \hline
 M_1 &\ \mZ/2\ &\ \mZ/4\ & \ \mZ/2\+\mZ/2\ & \ \mZ/2\+\mZ/4\ \\
 \hline
\ M_t\ (t>1)\ & \mZ/2 &\ \mZ/2\+\mZ/2\ &\ \mZ/4\+\mZ/2\ &\ \mZ/2\+\mZ/2\+\mZ/2\ \\  
 \hline 
 \end{array}  
 $$ 
 \end{table}
 It is convenient to organize this result in the form of Table~\ref{tab4} below, as in \cite{bd3}.

\begin{table}[!ht]
\caption{}\label{tab4} 

 $$ 
   \begin{array}{|c|ccccccc|ccccccc|}
\hline
  && ^1\* & ^2\* & ^3\* & \dots & ^\8\* &&& *^\8 & \dots & *^3 & *^2 & *^1 &\\
\hline
 \ \*_1 &&\ \mZ/2\ &\ \mZ/2\ &\ \mZ/2\ &\dots &\ \mZ/2\ &&& 0 & \dots & 0 & 0 & 0 &\\
 \ \*_2 &&\ \mZ/4\ & \mZ/2 & \mZ/2 & \dots & \mZ/2 &&&\ \mZ/2\ & \dots &\ \mZ/2\ &\ \mZ/2\ & 0 &\\
 \ \*_3 &&\ \mZ/4\ & \mZ/2 & \mZ/2 & \dots & \mZ/2 &&&\ \mZ/2\ & \dots &\ \mZ/2\ &\ \mZ/2\ & 0 &\\
 \ \vdots && \hdotsfor5 &&& \hdotsfor5 &\\
 \ \*_\8 &&\ \mZ/4\ & \mZ/2 & \mZ/2 & \dots & \mZ/2 &&&\ \mZ/2\ & \dots &\ \mZ/2\ &\ \mZ/2\ & 0 &\\
 \hline 
 \ _\8* && 0 & 0 & 0 & \dots & \ \mZ\ &&& \mZ/2 & \dots & \mZ/2 & \mZ/2 & \mZ/2 &\\
 \vdots && \hdotsfor5 &&& \hdotsfor5 &\\
 \ _3* && 0 & 0 & 0 & \dots & 0 &&&  \mZ/2 & \dots & \mZ/2 & \mZ/2 & \mZ/2 &\\
 \ _2* && 0 & 0 & 0 & \dots & 0 &&&  \mZ/2 & \dots & \mZ/2 & \mZ/2 & \mZ/2 &\\
 \ _1* && 0 & 0 & 0 & \dots & 0 &&&  \mZ/4 & \dots & \mZ/4 & \mZ/4 & \mZ/2 &\\
\hline
 \end{array}  
 $$ 
\end{table}
 In this table the row marked by $\*_t$ (respectively, $_t*$) shows the part of the group $\hos(\Si M_r,M_t)$ that comes from
 $\hos(\Si M_r,S^2)$ (respectively, from $\hos(\Si M_r,S^3)$\,). In the same way, the column marked by $^r\*$
 (respectively, $*^r$) shows the part of this group that comes from $\hos(S^3,M_t)$ (respectively, from $\hos(S^4,M_t)$\,).
 The columns $^\8\*$ and $*^\8$ correspond, respectively, to $\hos(S^4,\_\,)$ and $\hos(S^3,\_\,)$; the rows
 $\*_\8$ and $_\8*$ correspond, respectively, to $\hos(\_ \,,S^2)$ and $\hos(\_ \,,S^3)$.

 Therefore we consider the elements from $\El(\sW)$ as block matrices $(W^x_y)$,
 where $x\in\set{^r\*,*^r}$, $y\in\set{\*_t,^t*}$ and the block $W^x_y$ is with entries from the
 corresponding cell of Table~\ref{tab4}. Moreover, morphisms between Moore spaces induce the
 following transformations of vertical stripes $W^x$ and horizontal stripes $W_y$ of such a matrix,
 which we call \emph{admissible transformation}:
 \begin{itemize}
\item[$(\mathsf a)$] replacing the stripes $M^{^r\*}$ and $M^{*^r}$ by $M^{^r\*}X$ and $M^{*^r}X$;
 \item[$(\mathsf a')$]   replacing the stripes $M_{\*_t}$ and $M_{_t*}$ by $XM_{\*_t}$ and $XM_{_t*}$;
 \item[$(\mathsf b)$]  replacing $M^{*^r}$ by $M^{*^r}+M^{*^{r'}}X+M^{^s\*}Y$, where $r'>r$, $s$ arbitrary;
 \item[$(\mathsf b')$]  replacing $M_{\*_t}$ by $M_{\*_t}+XM_{\*{t'}}+YM_{_s*}$, where $t'>t$, $s$ arbitrary;
 \item[$(\mathsf c)$]  replacing $M^{^r\*}$ and $M^{*^r}$ by $M^{^r\*}+M^{^{r'}\*}X$ and
 $M^{*^r}+2^{r-r'}M^{*^{r'}}X$, where $r'<r$;
 \item[$(\mathsf c')$]  replacing $M_{_t*}$ and $M_{\*_t}$ by $M_{_t*}+XM_{_{t'}*}$ and
 $M_{\*_t}+2^{t-t'}XM_{\*_{t'}}$, where $t'<t$;
 \item[$(\mathsf d)$]  replacing $M^{^1\*}$ by $M^{^1\*}+2M^{^r\*}X+2M^{*^s}Y$; $r,s$ arbitrary;
 \item[$(\mathsf d')$]  replacing $M_{_1*}$ by $M_{_1*}+2XM_{_t*}+2YM_{\*_s}$; $r,s$ arbitrary;
 \item[$(\mathsf e)$]   replacing $M^{*^r}_{_1*}$ by $M^{*^r}_{_1*}+2M^{^s\*}_{\*_1}X$; $s$ arbitrary;
 \item[$(\mathsf e')$]  replacing $M^{^1\*}_{\*_t}$ by $M^{^1\*}_{\*_t}+2XM^{*^1}_{_s*}$; $s$ arbitrary.
\end{itemize}
 Here $X,Y$ denote arbitrary integer matrices of the appropriate size; in the
 transformations of types $(\mathsf a)$ and $(\mathsf a')$ the matrix $X$ must be invertible. 
 Two matrices $W,W'$  are isomorphic in $\El(\sW)$ \iff $W$ can be transformed to $W'$
 using admissible transformations.

 It is convenient first to reduce the block $W^{^\8\*}_{_\8*}$ to a diagonal form $D=\mathrm{diag}(\lst am)$
 with $a_1|a_2|\dots|a_m$. Let $a_k=2^{d_k}b_k$ with odd $b_k$. Denote by $W^{^{\8k}\*}$ and
 $W_{_{\8k}*}$ the parts of the stripes $W^{^\8\*}$ and $W_{_\8*}$ corresponding to the columns and rows
 with $d_k=d$ ($k=\8$ if $d_k=0$). Since all other matrices of these stripes are with entries from $\mZ/2$, we can
 make the parts $W^{^{\80}\*}$ and $W_{_{\80}*}$ zero. Moreover, using admissible transformations
 that do not change the block $D$, we can replace $W^{^{\8k}\*}$ by $W^{^{\8k}\*}+W^{^{\8l}\*}X$ and
 $W_{_{\8k}*}$ by $W_{_{\8k}*}+YW_{_{\8l}*}$ for any $l<k$. In what follows we always suppose that
 $W$ is already in this form.
 
 Call two matrices of this form $W,W'$ \emph{2-equivalent}, if there is a matrix $W''\iso W$ such that
 $W''\equiv W\mod2$. One can easily see that the problem of 2-equivalence of matrices from
 $\El(\sW)$ is actually a sort of \emph{bunch of chains} in the sense of \cite{bo,d3}. We use the
 paper \cite{d3} as the source for the further discussion. Namely, we have
 the chain $\kE=\set{\*_t,_t*,{_{\8k}*}}$ for the rows and the chain $\kF=\set{^r\*,*^r,{^{\8k}\*}}$ for the columns,
 where 
 \begin{multline*} 
   {_1*}<{_2*}<{_3*}<\dots<{_{\8\8}}<\dots<{_{\83}*}<{_{\82}*}<{_{\81}*}< \\
	  <\*_\8<\dots<\*_3<\*_2<\*_1,
 \end{multline*}
 \begin{multline*} 
  {^1\*}<{^2\*}<{^3\*}<\dots<{^{\8\8}\*}<\dots{^{\83}\*}<{^{\82}\*}<{^{\81}\*}<
	\\  <*^\8<\dots<*^3<*^2<*^1.
 \end{multline*} 
 The equivalence relation $\sim$ on $\kX=\kE\cup\kF$ is given by the rule
 $$ 
  \*_t\sim {_t*}\ (t\ne\8),\ ^r\*\sim *^r\ (r\ne\8),\  {^{\8k}\*}\sim{_{\8k}*}
 $$ 
 for all possible values of $t,r$ and $k\ne\8$. Thus we can get a classification of our matrices up to 2-equivalence
 from \cite{d3}. Namely, we write $x-y$ if either $x\in\kE,\,y\in\kF$ or vice versa, at least
one of them belongs to $\set{\*_t}\cup\set{*^r}$, moreover, $\set{x,y}\ne\set{\*_t,*^1}$ and
 $\set{x,y}\ne\set{\*_1,*^r}$. We call an \emph{$\kX$-word} a sequence
 $w=x_1\rho_2x_2\rho_3\dots \rho_nx_n$, where $x_i\in\kX,\,\rho_i\in\set{\sim,-}$, $\rho_i\ne \rho_{i+1}\ (i=2,\dots,n-1)$
 and $x_{i-1}\rho_ix_i$ holds in $\kX$ for all $i=2,\dots,n$. Such a word is called \emph{full} if the following
 conditions hold:
\begin{itemize}
\item either $\rho_2=\sim$ or $x_1\not\sim y$ for all $y\in\kX,\ y\ne x_1$;
 \item   either $\rho_n=\sim$ or $x_n\not\sim y$ for all $y\in\kX,\ y\ne x_n$.
\end{itemize}
 $w$ is called a \emph{cycle} if $\rho_2=\rho_n=-$ and $x_n\sim x_1$ in $\kX$. If, moreover, $w$ cannot be
 written in the form $v\sim v\sim \dots \sim v$ for a shorter word $v$, it is called \emph{aperiodic}. 
 We call a polynomial $f(t)\in\mZ/2[t]$ \emph{primitive} if it is a power of an irreducible polynomial
 with the leading coefficient $1$. We shall identify any word $w$ with its inverse and any cycle $w$
 with any of its cyclic shifts. Then the set of indecomposable representations of this bunch of
 chains is in 1-1 correspondence with the set $\kS\cup\kB$, where $\kS$ is the set of full words
 (up to inversion) and $\kB$ is the set of pairs $(w,f)$, where $w$ is an aperiodic cycle (up to a cyclic shift)
 and $f\ne t^d$ is a primitive polynomial. We call representations corresponding to $\kS$ \emph{strings}
 and those corresponding to $\kB$ \emph{bands}. 

 Note that an $\kX$-word can contain at most one element $^{\8k}\*$, at most one element $*_{\8k}$
 and at most one subword of the form $\*_t-*^r$ or its inverse. Replacing $w$ by its inverse,
 we shall suppose that there are no words of the form $*^r-\*_t$ or $^{\8k}\*\sim_{\8k}*$.
 It is convenient to rewrite this answer in a modified form. Namely, we replace the subword
 $_{\8k}*\sim ^{\8k}\*$, if it occurs, by $_k\eps^k$, also omit $x_1$ if $\rho_2=\sim$,
 omit $x_n$ if $\rho_n=\sim$ and omit all remaining symbols $\sim$.
 Then we replace every subword $^r\*-\*_t$ by $^r\*_t$, $\*_t-^r\*$ by $_t\*^r$, $_t*-*^r$
 by $_t*^r$, $*^r-_t*$ by $^r*_t$ and $\*_t-*^r$ by $_t\theta^r$.
 Note that in the last case $r\ne1$ and $t\ne1$. We also omit all signs $\sim$, replace any double
 superscript $^{rr}$ by $^r$ and any double subscript $_{tt}$ by $_t$. 
 Certainly, the original word can be easily restored from  such a shortened form.
 Now, any full word or its inverse can be written as a subword
 of one of the following words:
 \begin{align*}
&  ^{r_1}\*_{t_1}*^{r_2}\*_{t_2}*\dots ^{r_n}\*_{t_n} \quad \text{ (``usual word'')},  \\
&  _{t_{-m}}\*\dots ^{r_{-2}}*_{t_{-2}}\*^{r_{-1}}*_{t_{-1}}\theta^{r_1}\*_{t_1}*^{r_2}\*_{t_2}\dots *^{r_n}
 \quad\text{ (``theta-word'')}, \\
&  _{t_{-m}}\*\dots ^{r_{-2}}*_{t_{-2}}\*^{r_{-1}}*_k\eps^k\*_{t_1}*^{r_2}\*_{t_2}\dots *^{r_n}
 \quad\text{ (``epsilon-word'')},
 \end{align*} 
 Moreover, 
 \begin{itemize}
 \item  $\8$ can only occur at the ends of a word, not in a theta-word or epsilon-word.
  \item  In any theta-word $t_{-1}\ne 1$ and $r_1\ne1$.
\end{itemize}

 Any cycle or its shift can be written as 
 $$ 
  ^{r_1}\*_{t_1}*{r_2}\*_{t_2}*\dots ^{r_n}\*_{t_n}*^{r_1}.  
 $$ 
 
 The description of the representations in \cite{d3} also implies the following properties.

 \begin{proposition}\label{d71} 
 \begin{enumerate}
 \item  Any row (column) of a string contains at most $1$ non-zero element.
 \item  There are at most $2$ zero rows or columns in a string, namely, they are in the
 following stripes:
	\begin{enumerate}
	 \item  $M_{_t*}$ if $w$ has an end $\*_t$ (or $_t\*$), $t\ne\8$;
	 \item  $M^{*^r}$ if $w$ has an end $^r\*, r\ne\8$;
	 \item  $M_{\*_t}$ if $w$ has an end $_t*$;
	 \item  $M^{^r\*}$ if $w$ has an end $*^r$ (or $^r*$);
	 \item  $M_{_{\8k}*}$ if the left end of $w$ is $_k\eps^k$;
	 \item  $M^{^{\8k}\*}$ if the right end of $w$ is $_k\eps^k$.
	\end{enumerate}
 We call each end occurring in this list a \emph{distinguished end}.
 \item  The horizontal and vertical stripes of a band can be subdivided in such a way that every new 
 horizontal or vertical band has exactly $1$ non-zero block, which is invertible.
\end{enumerate}
   \end{proposition}

 Recall that elements modulo 4 only occur in the stripes $W^{^1\*}$ and $W_{_1\*}$.

  \begin{corollary}\label{d72} 
 Let $W\in\El(\sW)$ (with diagonal $W^{^\8\*}_{_\8*}$). Denote by $\ol W$ its reduction modulo 2 and by 
 $\ti W$ the matrix obtained from $W$ by replacing all invertible entries with $0$ (thus all entries of
 $\ti W$ are even). Suppose that $\ol W=\b+_{i=1}^m\ol W_i$,
 where all $\ol W_i$ are strings or bands. Then $W\iso W'$, where $\ol W'=\ol W$ and the only non-zero
 rows and columns of $\ti W'$ can be those corresponding to the distinguished ends of types 2(a-d)
 of Proposition~\ref{d71}. In particular, if some of $\ol W_i$ is a band, a theta-string or an epsilon-string,
 $W'$ has a direct summand $W_i$ such that $\ol W_i\equiv W_i\mod2$ and $\ti W_i=0$.
 \end{corollary}

 Thus we only have now to consider the case, when $W=W'$ and every $\ol W_i$ is a usual string.
 Suppose that $W_i$ corresponds to a string $w_i$. It is easy to verify that if $w_i$ and $w_j$ have a common
 distinguished end, there is a sequence of distinguished transformations, which does not change $\ol W$ and adds
 the row (or column) corresponding to this end in $\ti W_i$ to the row (or column) corresponding to this end in $\ti W_j$
 or vice versa. Hence, such rows (columns) are in some sense linearly ordered. As a consequence, we can transform
 $\ti W$ to a matrix having at most one non-zero element in every row and every column (without changing $\ol W$).
 It gives us the following description of indecomposable matrices from $\El(\sW)$ with $\ti W\ne0$.

  \begin{corollary}\label{d73} 
  Suppose that $W$ is an indecomposable matrix from $\El(\sW)$, such that $\ti W'\ne0$ for every matrix
 $W'\iso W$. Let $\ol W=\b+_{i=1}^m\ol W_i$, where each $\ol W_i$ is a usual string.
 There are, up to isomorphism, the following possibilities:
 \begin{enumerate}
\item  $m=1$, $\ol W$ corresponds to a word $w$ and $\ti W$ has a unique non-zero element in the block
 $W^a_b$ for the following choices: 
 \begin{align*} 
 &	w={_{t_1}*}^{r_2}\*_{t_2}\dots,\ a={^1\*},\ b=\*_{t_1}\ (t_1=\ne 1); \tag{a}\\
 &	w={^{r_1}\*}_{t_1}*^{r_2}\*\dots, \ a=*^{r_1},\ b={_1*}\ (r_1\ne1); \tag{b},\\
 &	w={_1*}^{r_1}\*_{t_1}*^{r_2}\dots,\ a=*^r,\ b={_1*}; \tag{c},\\
 &	w={^1\*}_{t_1}*^{r_2}\*_{t_2}\dots,\ a={^1\*},\ b=\*_t; \tag{d}.
 \end{align*} 
 \item  $m=2$, $\ol W_i\ (i=1,2)$ correspond to the words $w_i$ and $\ti W$ has a unique non-zero element
 in the block $W^a_b$, where 
 \begin{align*} 
 &	w_1={_1*}^{r_{-1}}\*_{t_{-1}}*^{r_{-2}}\*\dots,\ w_2={^{r_1}\*}_{t_1}*^{r_2}\*_{t_2}\dots,\ 
	a=*^{r_1},\  	b={_1*},  \tag{e}\\
 &	w_1={^1\*}_{t_1}*^{r_1}\*_{t_2}*\dots,\ w_2={_{t_{-1}}}*^{r_{-1}}\*_{t_{-2}}*^{r_{-2}}\*\dots,\ 
	a={^1\*},\ b=\*_{t_{-1}}. \tag{f}
 \end{align*} 
\end{enumerate}
 We encode these matrices by the following words $w$:
 \begin{align*} 
  &  w=\dots \*^{r_2}*_{t_2}\*^{r_1}*_{t_1}\theta^1 \quad\text{in case (a)},\\
  &  w={_1\theta}^{r_1}\*_{t_1}*^{r_2}\*_{t_2}*\dots \quad\text{in case (b)},\\
  &  w=\dots *_{t_2}\*^{r_2}*_{t_1}\*^{r_1}*_1\theta^r \quad\text{in case (c)},\\
  &  w={_t\theta}^1\*_{t_1}*^{r_2}*_{t_2}\*^{r_3}*\dots \quad\text{in case (d)},\\
  &  w=\dots *^{r_{-2}}\*_{t_{-2}}*^{r_{-1}}*_1\theta^{r_1}\*_{t_1}*^{r_2}\*\dots \quad\text{in case (e)},\\
  &  w=\dots *{t_{-2}}\*^{r_{-1}}*_{t_{-1}}\theta^1\*_{t_1}*^{r_2}\*_{t_2}*\dots \quad\text{in case (f)},\\
 \end{align*} 
 We call these words ``theta-words'' as well.
 \end{corollary}

 Obviously, cases (a-d) always give indecomposable matrices. On the other hand, one can check that
 in case (e) $W$ is indecomposable \iff  $(r_{-1}+1,t_1,r_{-2},t_2,\dots)<(r_1,t_{-2},r_2,t_{-3},\dots)$
 with respect to the lexicographical order \cite{bd3}. In case (f) $W$ is indecomposable \iff 
 $(t_1+1,r_{-1},t_2,r_{-2},\dots)<(t_{-1},r_2,t_{-2},r_3,\dots)$ lexicographically. Thus we obtain a complete
 list of non-isomorphic indecomposable matrices from $\El(\sW)$. Moreover, it is easy to verify that
 they remain pairwise non-isomorphic and indecomposable in $\El(\sW)\kJ$ as well. Thus,
 using Theorem~\ref{bh-1}, we get the following result.

  \begin{theorem}[Baues--Hennes \cite{bh}]\label{d74} 
  Indecomposable polyhedra from $\cw^3_4$ are in 1-1 correspondence with usual words, theta-words,
 epsilon-words and bands defined above, with the only restriction that in a theta-word
 $w=\dots ^{r_{-2}}*_{t_{-2}}\*^{r_{-1}}*_{t_{-1}}\theta^{r_1}\*_{t_1}*^{r_2}\*_{t_2}\dots$ the following
 conditions hold: \em
 \begin{align*} 
 \text{ if $\,t_{-1}=1$, then }\ &(r_{-1}+1,t_1,r_{-2},t_2,\dots)<(r_1,t_{-2},r_2,t_{-3},\dots), \\
 \text{ if $\,r_1=1$, then }\ &(t_1+1,r_{-1},t_2,r_{-2},\dots)<(t_{-1},r_2,t_{-2},r_3,\dots)
 \end{align*} 
 (lexicographically in both cases).
 \end{theorem}

 The gluings of spheres corresponding to these words can be described as follows:
 $$
  \xymatrix@R=.1ex@C=2ex@!R{
	 {7} \ar@{.} [rrrrrrrrrrrrrr] &&&& \cel \ar@{-}[dd] &&&& \cel \ar@{-}[dd]
 &&&& \cel \ar@{-}[dd] && {} \\
  {}\\
	{6} \ar@{.} [rrrrrrrrrrrrrr] &&&& \cel \ar@{-}[ddddrr] &&&& \cel \ar@{-}[ddddrr] 
 &&&& \cel \ar@{-}[ddr] && {}  \\
	& {\cdots\quad } \\
	{5} \ar@{.} [rrrrrrrrrrrrrr] & *=0{}& \cel \ar@{-}[uuuurr] \ar@{-}[dd] &&&&
 \cel \ar@{-}[uuuurr] \ar@{-}[dd] &&&&  \cel \ar@{-}[uuuurr] \ar@{-}[dd] &&& *=0{} & {} \\
	&&&&&&&&&&&&& {\quad \cdots}  \\
	{4}  \ar@{.} [rrrrrrrrrrrrrr] && \cel@\ar@{-}[uul] &&&& \cel &&&& \cel &&&& {}
 }
 $$
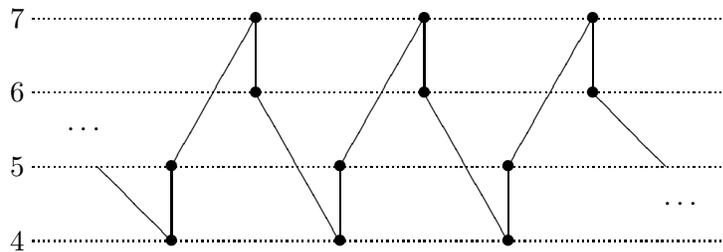
\centerline{\bf for a usual word}

\medskip
 $$
  \xymatrix@R=.1ex@C=2ex@!R{
	 {7} \ar@{.} [rrrrrrrrrrrrrrrr] && \cel \ar@{-}[dd] \ar@{-}[ddddrr] &&&&
 \cel \ar@{-}[dd] \ar@{-}[ddddrr] &&&& \cel \ar@{-}[dd] &&&& \cel \ar@{-}[dd] && {} \\
 {} \\
	{6} \ar@{.} [rrrrrrrrrrrrrrrr] && \cel \ar@{-}[ddl] &&&& \cel \ar@{-}[ddddll] 
 &&&& \cel \ar@{-}[ddddrr] &&&& \cel \ar@{-}[ddr] && {}  \\
	{} \\
	{5} \ar@{.} [rrrrrrrrrrrrrrrr] & *=0{}&&& \cel \ar@{-}[dd] &&&&
 \cel \ar@{-}[dd] &&&&  \cel \ar@{-}[uuuurr] \ar@{-}[dd] &&& *=0{} & {} \\
	& {\cdots\quad} &&&&&&&&&&&&&& {\quad \cdots}  \\
	{4}  \ar@{.} [rrrrrrrrrrrrrrrr] &&&& \cel &&&& \cel \ar@{-}[uuuuuurr] &&&& \cel &&&& {}
 }
 $$
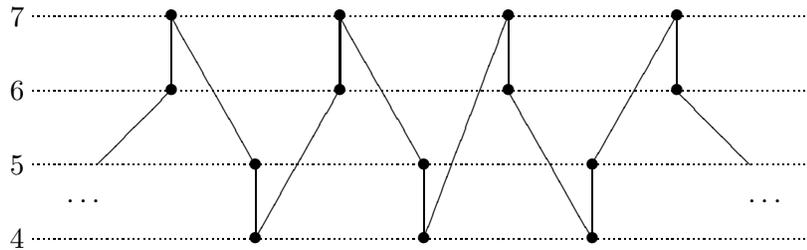
\centerline{\bf for a  theta-word}

\medskip
 $$ 
  \xymatrix@R=.1ex@C=2ex@!R{
	 {7} \ar@{.} [rrrrrrrrrrrrrrrr] &&&& \cel \ar@{-}[dd] \ar@{-}[ddddrr] &&&&&&
  \cel \ar@{-}[dd] \ar@{-}[ddddll] &&&& \cel \ar@{-}[dd] \ar@{-}[ddddll]  & & {} \\
  &{\cdots\quad} \\
	{6} \ar@{.} [rrrrrrrrrrrrrrrr] &&&& \cel \ar@{-}[ddddll] &&&&
 \cel \ar@{-}[ddddll] \ar@{-}[dd] &&  \cel \ar@{-}[ddddrr]  &&&& \cel \ar@{-}[ddr] && {} \\
  {} \\
	{5} \ar@{.} [rrrrrrrrrrrrrrrr] && \cel \ar@{-}[uul] \ar@{-}[dd] &&&& \cel \ar@{-}[dd] 
 && \cel  &&&& \cel \ar@{-}[dd] &&& *=0{}& {}  \\
	& &&&&&&&&&&&&&& {\quad \cdots}  \\
	{4}  \ar@{.} [rrrrrrrrrrrrrrrr] && \cel &&&& \cel  &&&&&& \cel &&&& {}
 }    
 $$ 
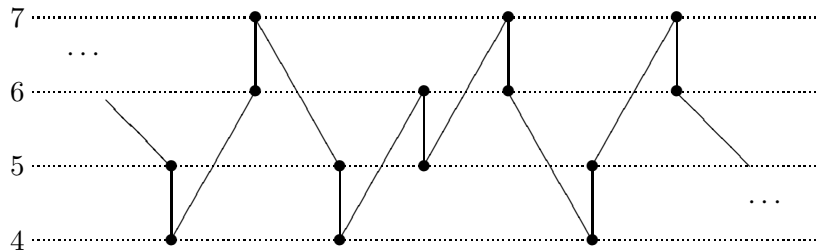
\centerline{\bf for an epsilon-word}

\medskip\noindent
 In these diagrams vertical segments present the suspended atoms $M_r$, slanted lines correspond to
 the gluings arising from Hopf maps $S^{d+1}\to S^d$, while the long slanted line in a theta-word
 shows the gluing arising from the doubled Hopf map $S^6\to S^4$. 

 Note that all atoms from $\cw^3_4$ are $p$-primary ($2$-primary, except $M(q)$ with odd $q$).
 Therefore, we have the uniqueness of decomposition of spaces from $\cw^3$ into bouquets
 of suspended atoms.

 \section{Bigger dimensions. Wildness}
 \label{bdw} 

 Unfortunately, if we pass to bigger dimensions, the calculations as above become extremely complicated.
 In the representations theory the arising problems are usually called ``\emph{wild}.'' 
 Non-formally it means that the classification problem for a given category
 contains the classification of representations of arbitrary (finitely generated) algebras over a field.
 It is well-known, since at least 1969 \cite{gp}, that it is enough to show that this problem contains
 the classification of pairs of linear mappings (up to simultaneous conjugacy), or, equivalently,
 the classification of triples of linear mappings
 \begin{equation}\label{bde1} 
  \xymatrix{{V_1\,} \ar[r] \ar@<1ex>[r] \ar@<-1ex>[r] & {\ V_2}}
 \end{equation}
 On the other hand, the problems
 like one considered in the preceding section, where indecomposable objects can be parameterised by
 several ``discrete,'' or combinatorial parameters (as $\kX$-words above) and at most one ``continuous''
 parameter (as a primitive polynomial in the description of bands), are called ``\emph{tame}.''
 The problems, where the answer is purely combinatorial, like the classification of atoms of
 dimensions $d\le 5$, are called ``finite.'' I shall not precise these notions formally. The reader can
 consult, for instance, the survey \cite{d2}, where it is done within the framework of representation
 theory. An important question in the representation theory is to distinguish finite, tame
 and wild cases. The following result accomplishes such an investigation for stable homotopy types.

  \begin{proposition}\label{bdw1} 
  The classification problem for the category $\cw^4$ is wild.
 \end{proposition}  
 \begin{proof} 
  Let $\bB$ be the category of bouquets of Moore atoms $M=M_1$, $\bA=\Si^2\bB$. Then $\cw^4$
 contains the subcategory $\Si^3(\bA\dagg\bB)\iso\bA\dagg\bB$. Corollary~\ref{bh-2} shows that the
 category $\bA\dagg\bB$ is  representation equivalent to $\El(\sH)$, where $\sH$ is the restriction
 of $\hos$ onto $\bA^\circ\xx\bB$. We know that $\hos(M,M)=\mZ/4$. Therefore, we only have to
 show that $\hos(\Si^2M,M)\iso\mZ/2\mZ/2\xx\mZ/2$. Indeed, it implies the category $\El(\sH)$ is
 representation equivalent to the category of diagrams of the shape \eqref{bde1}. 

 The cofibration sequence $S^2\str2S^2\to M\to S^3\str2 S^3$ and the Hopf map $\eta:S^5\to S^4$
 produce the following commutative diagram:
 $$ 
   \begin{CD}
      0 @>>> \mZ/2 @>>> \pi^S_4(M) @>>> \mZ/2 @>>> 0 \\
   &&	@V\eta^*VV @VVV  @VV\wr V \\
     0 @>>> \mZ/2 @>>> \pi^S_5(M) @>>> \mZ/2 @>>> 0 ,
 \end{CD}  
 $$ 
 Since $\eta^3=4\nu$, where $\nu$ is the element of order $8$ in $\pi^S_5(S^2)=\mZ/24$ \cite{tod},
 actually $\eta^*=0$, so the lower row splits and $\pi^S_5(M)=\mZ/2\+\mZ/2$. Just in the same way
 we show that $\hos(\Si^2M,S^2)=\mZ/2\+\mZ/2$. Now, applying the functors
 $\hos(\_\,,S^2)$ and $\hos(\_\,,M)$ to the same cofibration sequence, we get the commutative diagram
 $$ 
   \begin{CD}
         0 @>>> \mZ/2 @>>> \hos(\Si^2M,S^2) @>>> \mZ/2 @>>> 0 \\
   &&	@VVV @VVV  @VV\wr V \\
     0 @>>> \mZ/2\+\mZ/2 @>>> \hos(\Si^2M,M) @>>> \mZ/2 @>>> 0 .
 \end{CD}  
 $$ 
 Since the upper row of this diagram splits, the lower one splits as well, hence
 $\hos(\Si^2M,M)=\mZ/2\+\mZ/2\+\mZ/2$. It accomplishes the proof.
 \end{proof} 

 We can summarize the obtained results in the following theorem.

 \begin{theorem}
    The category $\cw^k$ is of finite type for $k\le 2$, tame for $k=3$ and wild for $k\ge4$.
 \end{theorem}

 \section{Torsion free atoms. Dimension 9.}
 \label{tf} 

 Nevertheless, if we consider \emph{torsion free} atoms, the situation becomes much simpler. Namely, in this case
 neither sphere of dimension $d$ can be attached to the spheres of dimension $d-1$, thus in the picture describing
 the gluing of spheres there is no fragments of the sort 
 $$ 
 \xymatrix@=1ex{ d &{} \ar@{.}[rr] & \cel\ar@{-}[d] &{}\\ {d-1} &{}\ar@{.}[rr] &\cel &{} }
 $$ 
 Therefore, a calculation of atoms from $\cwf^k_{k+1}$ can be organized as follows. Denote by
 $\bB_k$ the full subcategory of $\cw$  consisting of bouquets of torsion free suspended atoms 
 of dimension $2k$ and by $\bS_k$ the category of bouquets of spheres $S^{2k}$. Let $\Ga_m(X)$
 denote the subgroup $\im\{\pi_m^S(X^{m-1})\to\pi_m(X)\}$ of $\pi_m^S(X)$. When $X$ runs
 through $\bB_k$, $\Ga_{2k}$ can be considered as an $\bS_k$-$\bB_k$-bimodule; we denote this
 bimodule by $\sG_k$. Then the following analogue of Theorem~\ref{bh-1} holds (with essentially
 the same proof).

  \begin{proposition}\label{tf-1} 
  Denote by $\kI$ the ideal of the category $\cwf^k_{k+1}$ consisting of all morphisms $X\to X'$
 that factor through an object from $\bB_k$, and by $\kJ$ the ideal of the category $\El(\sG_k)$
 consisting of such morphisms $(\al,\be):f\to f'$ that $\al$ factors through $f$ and $\be$ factors
 through $f'$. Then $\cwf^k_{k+1}/kI\iso \El(\Ga_k)/kJ$. Moreover, both $\kI^2=0$ and
 $\kJ^2=0$, hence the categories $\cwf^k_{k+1}$ and $\El(\sG_k)$ are representation
 equivalent.
 \end{proposition} 
 \begin{proof} 
 The only new claim here is that $\kJ^2=0$. But this equality immediately follows from the fact
 that if a morphism $X\to S^m$ factors through $X^{m-1}$, it is zero.
 \end{proof} 

 Thus a torsion free atom of dimension 7 can be obtained as a cone of a map $f:mS^6\to Y$,  where $Y$ is
 a bouquet of spheres $S^4,S^5$ and suspended Chang atoms $C^6(\eta)$, while $f\in\Ga_6(Y)$. Easy
 calculations, like above, give the following values of $\Ga_6$:

\medskip
 $$
 \begin{array}{|c|c|c|c|}
  \hline
  X &  S^4 &  S^5 & \ C^6(\eta)\ \\
 \hline
 \ \Ga_6\ &\ \mZ/2\ &\ \mZ/2\ & 0 \\
 \hline
 \end{array} 
 $$

\medskip\noindent
 (The last $0$ is due to the fact that the map $\eta_*:\pi_6(S^5)\to\pi_6(S^4)$ is an epimorphism \cite{tod}).
 The Hopf map $\eta:S^5\to S^4$ induces an isomorphism $\Ga_6(S^5)\to\Ga_6(S_4)$. Therefore, the only indecomposable
 torsion free atom of dimension 7 is the gluing $C(\eta^2)=C^7(\eta^2)=S^4\cup_{\eta^2}B^7$.
 (Note that such an atom must contain at least one $4$-dimensional cell.) Moreover, all torsion free atoms of
 dimensions $d\le7$ are $2$-primary.

\bigskip
 A torsion free atom of dimension 9 is a cone of some map $f:mS^8\to Y$, where $Y$ is a bouquet
 of spheres $S^i\ (5\le i\le 7)$, suspended Chang atoms $C^7(\eta),C^8(\eta)$ and suspended atoms $C^8(\eta^2)$. One can calculate
 the following table of the groups $\Ga_8$ for these spaces:

\medskip
 $$ 
  \begin{array}{|c|c|c|c|c|c|c|}
  \hline
 X & S^5 & S^6 & S^7 &\ C^7(\eta)\ &\ C^8(\eta)\ &\ C^8(\eta^2) \ \\
 \hline
 \ \Ga_8\ &\ \mZ/24\ &\ \mZ/2\ &\ \mZ/2\ & \mZ/12 & 0 & \mZ/12 \\
 \hline
 \end{array} 
 $$ 

\medskip\noindent
 Morphisms between these spaces induce epimorphims $\Ga_8(S^5)\to\Ga_8(C^7(\eta^2))\to\Ga_8(C^7(\eta))$,
 $\Ga_8(S^7)\to\Ga_8(S_6)$ and monomorphisms $\Ga_8(S^7)\to\Ga_8(C^7(\eta)\to\Ga_8(S^5)$,
 $\Ga_8(S^6)\to\Ga_8(S^5)$. It can be deduced either from \cite{tod} or, perhaps easier, from the results
 of \cite{un}, cf. \cite{bd1}. (The only non-trivial one is the monomorphism $\Ga_8(S^7)\to\Ga_8(\Si^2A(\eta))\,$.)
 Again we consider the map $f$ as a block matrix 
 $$ 
  F=  \begin{pmatrix}
  F_1 & F_2 &\ F_3 &\ F_4 &\ F_6 
 \end{pmatrix}^\top .
 $$ 
 Here $F_i$ is of size $m_i\xx m$ with entries from $\Ga_8(Y_i)$, where  
 $$ 
 Y_i=    \begin{cases}
 \ S^5 &\text{if } i=1,\\ \ C^7(\eta) &\text{if } i=2, \\ \ C^8(\eta^2)&\text{if } i=3,\\
 \ S^6 &\text{if } i=4,\\ \ S^5 &\text{if } i=6.
 \end{cases} 
 $$ 
 We have written $F_6$, not $F_5$, in order to match the notations of the Example~\ref{bc-1a}; so we set $I_1=\set{1,2,3,4,6}$.
 Using the automorphisms of $mS^7$ and of $Y$, one can replace the matrix $F$ by $PFQ$, where $P\in\GL(m,\mZ)$ and
 $Q=(Q_{ij})_{i,j\in I_1}$ is an invertible integer block matrix, where the block $Q_{ij}$ is of size $m_i\xx m_j$ with the following
 restrictions for the entries $a\in Q_{ij}$:
 \begin{align*}
  a=&0 \quad\text{ for } i\in\set{4,6},\ j<i, \\
  a\equiv& 0 \mod2\, \text{ for } (ij)\in\set{(12),(13),(23)},\\
  a\equiv& 0 \mod6\, \text{ for } (ij)=(26),\\
  a\equiv& 0 \mod{12}\, \text{ for } j\in{4,6},\ i\in\set{1,2,3},\ (ij)\ne(26).
 \end{align*} 
 Thus we have come to the \emph{bimodule category} $\El(\sU_1)$ considered in Example~\ref{bc-1b}, so we can use
 Corollary~\ref{bc-3}, which describes all indecomposable objects of this category. Certainly, we are not interested
 in the ``empty'' objects $\0_i$, since they correspond to the spaces with no 9-dimensional cells.
  Note also that the matrices $(1_4),(1_6)$ correspond not to
 atoms, but to suspended atoms $C^9(\eta^2)$ and $C^9(\eta)$. 
 We use the following notation for the atoms corresponding to other indecomposable matrices $F$:
 $$ \def\for{&& \text{for} &&}
    \begin{array}{ccccc}
  A(v) \for  (v_1),\\
  A(\eta v)  \for (v_2),\\ \vspace*{.5ex}
  A(\eta^2 v) \for (v_3),\\ \vspace*{1ex}
  A(v\eta) \for \ds\binom{v_1}{1_6},\\ \vspace*{1ex}
  A(v\eta^2) \for \ds\binom{v_1}{1_4},\\  \vspace*{1ex}
  A(\eta v\eta) \for \ds\binom{v_2}{1_6},\\  \vspace*{1ex}
  A(\eta v\eta^2) \for \ds\binom{v_2}{1_4},\\  \vspace*{1ex}
  A(\eta^2 v\eta) \for \ds\binom{v_3}{1_6},\\  
  A(\eta^2 v\eta^2) \for \ds\binom{v_3}{1_4}.
 \end{array} 
 $$ 

 So we have proved 
 \begin{theorem}\label{tf-2}
  Every torsion free atom of dimension $9$ is isomorphic to one of the atoms $A(w)$ with
 $w\in\set{v,\eta v,\eta^2 v,v\eta,v\eta^2,\eta v\eta,\eta v\eta^2,\eta^2 v\eta,\eta^2 v\eta^2}$.
 \end{theorem} 

 Using the gluing diagrams, these atoms can be described as in Table~\ref{tab9} below.  
 \begin{table}[!ht]
 \caption{}\label{tab9}
 $$ 
   \xymatrix@R=.7ex@C=2ex{
  9  \ar@{.}[rrrrrrrrrr] & \cel \car[dddd] && \cel \car[ddddl] && \cel \car[ddddl] &&
	\cel \car[dd] \car[ddddl] && \cel \car[ddd] \car[ddddl] & {} \\
  8  \ar@{.}[rrrrrrrrrr] & &&& \cel \car[ddd] &&&&&& {} \\
  7  \ar@{.}[rrrrrrrrrr] & & \cel \car[dd] &&&&& \cel &&& {} \\
  6  \ar@{.}[rrrrrrrrrr] & &&&&&&&& \cel & {} \\
  5  \ar@{.}[rrrrrrrrrr] & \cel & \cel && \cel && \cel && \cel && {}	\\
  & A(v) & A(\eta v) && A(\eta^2 v) && A(v\eta) &&  A(v\eta^2) }
 $$ 

 $$ 
   \xymatrix@R=.7ex@C=2ex{ 
  9  \ar@{.}[rrrrrrrrrrrr] && \cel \car[ddddl]\car[dd] &&& \cel \car[ddddl]\car[ddd] &&& \cel \car[ddddl]\car[dd]
 &&&\cel \car[ddd] \car[ddddl]  & {} \\
  8  \ar@{.}[rrrrrrrrrrrr] & \cel \car[ddd] &&& \cel \car[ddd] &&&&&&&& {} \\
  7  \ar@{.}[rrrrrrrrrrrr] & & \cel &&&&& \cel \car[dd] &\cel && \cel\car[dd] &&{} \\
  6  \ar@{.}[rrrrrrrrrrrr] & &&&& \cel &&&&&& \cel & {} \\
  5  \ar@{.}[rrrrrrrrrrrr]  &\cel &&& \cel &&& \cel &&& \cel && {}	\\
  & A(\eta^2 v\eta) &&& A(\eta^2 v\eta^2) &&& A(\eta v\eta) &&&  A(\eta v\eta^2) }
 $$ 
 \end{table}

 One can also check that the $2$-primary atoms in this list are those with $v$ divisible by $3$, while the only $3$-primary one is
 $A(8)$. Thus there are altogether $29$ primary suspended atoms of dimension at most 9. The congruent ones are only
 $A(3)$ and $A(9)$. Indeed, $A(3)\vee S^5$ corresponds to the matrix $\ds \binom 30 \mod{24}$. But the latter can be easily transformed
 to $\ds\binom 90 \mod{24}$, which corresponds to $A(9)\vee S^5$:
 $$ 
  \binom 30 \to \binom 3{12} \to \binom {-9}{12} \to \binom 9{12} \to \binom 90.  
 $$ 
 (At the last step we add the first row multiplied by $4$ to the second one; all other transformations are obvious.) One can verify
 that all  other $2$-primary atoms are pairwise non-congruent. 
 
 \begin{corollary}\label{tf-3}
  The Grothendieck group $K_0(\cwf^4)$ is a free abelian group of rank 29.
 \end{corollary} 

 Note that the matrix presentations allows easily to find the images in $K_0(\cwf^4)$ of all atoms. For instance,
 the equivalence of matrices 
 $$ 
   \begin{pmatrix}  8_1 &0 \\ 0 & 3_1 \end{pmatrix}  \sim
   \begin{pmatrix}  1_1 &0 \\ 0 & 0_1 \end{pmatrix} 
 $$ 
 implies that $A(8)\vee A(3)\iso A(1)\vee S^5\vee S^9$, thus in $K_0(\cwf^4)$ we have
 $$ 
  [A(1)]=[A(8)]+[A(3)]-[S^5]-[S^9].  
 $$ 
 The reader can easily make analogous calculations for all atoms of Table~\ref{tab9}.

 \section{Torsion free atoms. Dimension 11.}
 \label{tf1} 

 For torsion free atoms of dimension 11 analogous calculations have been done in \cite{bd4}.
 Nevertheless, they are a bit cumbersome, so we propose here another, though
 rather similar, approach. Namely, denote by $\bS'_k$ the category of bouquets of spheres $S^{2k-1}$
 and $S^{2k}$, by $\bB'_k$ the category of bouquets of suspended atoms of dimension $2k-1$
 and by $\sG'_k$ the $\bS'_k$-$\bB'_k$-bimodule such that
 $$
 \sG'_k(S^{2k-1},B)=\Ga_{2k-1}(B)\quad\text{and}\quad
 \sG'_k(S^{2k},B)=\Ga_{2k}(B).
 $$
 
  \begin{proposition}\label{tf1-1} 
    Denote by $\kI'$ the ideal of the category $\cwf^k_{k+1}$ consisting of all morphisms $X\to X'$
 that factors through an object from $\bB'_k$, and by $\kJ'$ the ideal of the category $\El(\sG'_k)$
 consisting of such morphisms $(\al,\be):f\to f'$ that $\al$ factors through $f$ and $\be$ factors
 through $f'$. Then $\cwf^k_{k+1}/\kI'\iso \El(\Ga'_k)/\kJ'$. Moreover, both $(\kI')^2=0$ and
 $(\kJ')^2=0$, hence the categories $\cwf^k_{k+1}$ and $\El(\Ga'_k)$ are representation
 equivalent.
 \end{proposition}
 
 Thus we obtain torsion free atoms if dimension $11$ as cones of maps $S\to Y$, where $S$ is a bouquet
 of spheres of dimensions $9$ and $10$, while $Y$ is a bouquet of $5$-connected suspended atoms of
 dimensions $6\le d\le 9$. Note that at least one of these atoms must have a cell of dimension $6$
 in order that such a cone be an atom.

 Just as above, we have the following values of $\Ga_9$ and $\Ga_{10}$ for such atoms: 

 $$ 
    \begin{array}{|c|c|c|c|c|c|c|c|}
  \hline
 X & S^6 &\ C^8(\eta)\ &\ C^9(\eta^2) \ & S^7 &\ C^9(\eta)\ &  S^8 & S^9 \\
 \hline
 \ \Ga_9\ &\ \mZ/24\ &\ \mZ/12\ &\ \mZ/12\ &\ \mZ/2\ & 0&\ \mZ/2\ & 0 \\
 \hline
 \ \Ga_{10} & 0 &0&0 & \mZ/24 &\ \mZ/12\ &\ \mZ/2\ & \ \mZ/2\ \\
 \hline  
 \end{array} 
 $$ 

\medskip
\noindent
 (We have arranged this table taking into account the known maps between these groups, as above.)
 The Hopf map $S^{10}\to S^9$ induces monomorphisms in the 4th and the 6th columns of this table, while the
 maps between suspended atoms induce homomorphisms analogous to those of the preceding section. Thus a morphism
 $f:S\to Y$ can be described by a matrix 
 $$ 
  F=  \begin{pmatrix}
  F_1 & F_2 & F_3 & F_4 & 0 & F_6 & 0 \\ 0 & 0 & 0 & G_4 & G_5&G_6& G_7 
 \end{pmatrix}^\top ,
 $$ 
 where the matrix $F_i$ ($G_i)$ has entries from the first row (respectively, second row) and the $i$-th column of
 the table above. Two matrices, $F$ and $F'$, define homotopic polyhedra if $F'=PFQ$, where $P,Q$ are
 matrices over the tiled orders, respectively, 
 $$
 \begin{pmatrix}
  \ 1 & 2 & 2 & 12 & 24 & 12 & 24 \,\\ \ 1 & 1 & 1 & 12 & 24 & 6 & 24 \,\\ \ 1 & 2 & 1 & 12 & 24 & 12 & 24\, \\
\  0&0&0& 1 & 2 & 12^* & 12\, \\ \ 0&0&0& 1 & 1 & 12 & 6\, \\\ 0&0&0&0&0& 1&1\\ \ 0&0&0&0&0&0&1 \,
 \end{pmatrix} 
 \quad\text{and}\quad  
\begin{pmatrix} \ 1 & 12^*\, \\ \ 0&1\, \end{pmatrix}. \\
 $$
 Here $12^*$ shows that the corresponding element obeys the $\ast$-rule \eqref{e21}, i.e.
 induces a \emph{non-zero} map $\mZ/2\to\mZ/2$ and acts as usual multiplication by $12$ in all other cases.

 Thus we have obtained the \emph{bimodule category} $\El(\sU_2)$ from Example~\ref{bc-1a}, so can use the
 list of indecomposable objects from Theorem~\ref{bc-2}. Moreover, we only have to consider the matrices
 having non-empty $G$-column and one of the parts $F_1,F_2,F_3$ (otherwise we have no
 $11$-dimensional or no $6$-dimensional cells). Therefore, a complete list of atoms arises from the
 following matrices:
 $$ 
   \begin{pmatrix} v_i & 0 \\ 1_4 & w_4 \end{pmatrix} ,\quad
   \begin{pmatrix} v_i & 0 \\ 1_4 & w_4 \\ 0 & 1_6 \end{pmatrix} ,\quad
   \begin{pmatrix} v_i & 0 \\ 1_4 & w_4 \\ 0 & 1_7 \end{pmatrix} , 
 $$ 
 where $i\in\set{1,2,3}$, $v,w\in\set{1,2,3,4,5,6}$. We omit the upper indices of Theorem~\ref{bc-2},
 since here they coincide with the column number; the lower indices show to which horizontal stripe of the matrix $F$
 the corresponding elements  belong.
 It gives the following list of $11$-dimensional torsion free atoms.

 \begin{theorem}\label{tf1-2} 
  Every torsion free atom of dimension 11 is isomorphic to one of the atoms of Table~\ref{tab11} below.
 \end{theorem}

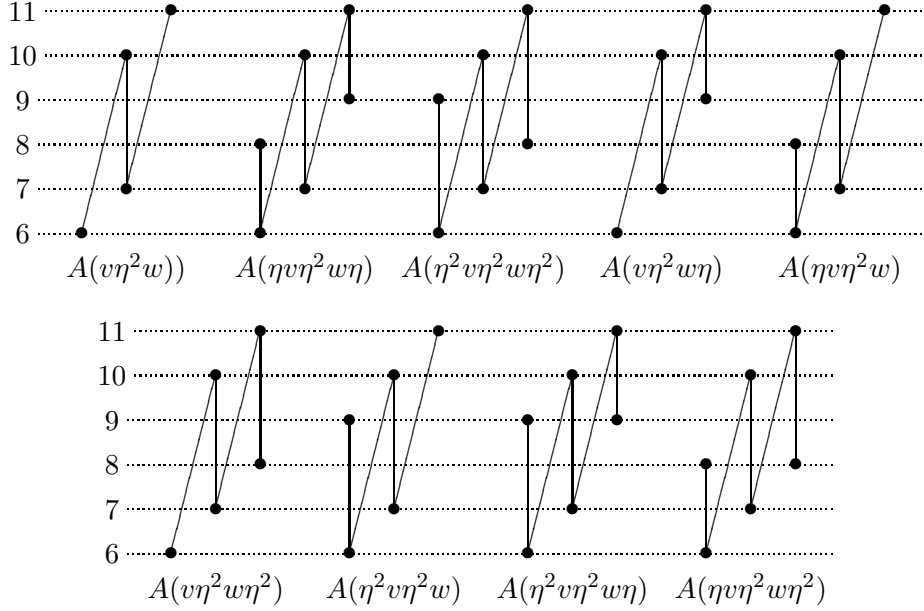
\begin{table}[!ht]
 \caption{}\label{tab11}
 $$ 
   \xymatrix@R=.7ex@C=2.2ex{
   11 \ar@{.}[rrrrrrrrrrrrrrrrrrrr] &&&\cel\ar@{-}[ddddl] &&&&\cel\ar@{-}[ddddl]\ar@{-}[dd]  &&&&
 \cel\ar@{-}[ddd] \ar@{-}[ddddl]&&&& \cel\ar@{-}[ddddl]\ar@{-}[dd]  &&&& \cel\ar@{-}[ddddl]  &{} \\
   10 \ar@{.}[rrrrrrrrrrrrrrrrrrrr] &&\cel\ar@{-}[ddddl]\ar@{-}[ddd]  &&&&\cel\ar@{-}[ddddl]\ar@{-}[ddd]
  &&&&\cel\ar@{-}[ddddl] \ar@{-}[ddd] &&&& \cel\ar@{-}[ddddl]\ar@{-}[ddd]  &&&& \cel\ar@{-}[ddd] \ar@{-}[ddddl]&&{} \\
   9 \ar@{.}[rrrrrrrrrrrrrrrrrrrr] &&&&&&&\cel &&\cel\ar@{-}[ddd] &&&&&&\cel &&&& &{} \\
   8 \ar@{.}[rrrrrrrrrrrrrrrrrrrr] &&&&&\cel\ar@{-}[dd] &&&&&&\cel &&&&&& \cel\ar@{-}[dd] &&& {} \\
   7 \ar@{.}[rrrrrrrrrrrrrrrrrrrr] &&\cel &&&&\cel &&&&\cel &&&& \cel &&&& \cel &&{} \\
   6 \ar@{.}[rrrrrrrrrrrrrrrrrrrr] & \cel &&&&\cel &&&&\cel &&&& \cel &&&& \cel &&&{} \\
   && *=0{A(v\eta^2w))} &&&& *=0{A(\eta v\eta^2w\eta)} &&&& *=0{A(\eta^2v\eta^2w\eta^2)}
 &&&& *=0{A(v\eta^2w\eta)} &&&& *=0{A(\eta v\eta^2w)} &&
 }
 $$

 $$ 
   \xymatrix@R=.7ex@C=2.2ex{
   11 \ar@{.}[rrrrrrrrrrrrrrrr] &&&\cel\ar@{-}[ddddl]\ar@{-}[ddd] &&&&\cel\ar@{-}[ddddl]  &&&&
 \cel\ar@{-}[dd] \ar@{-}[ddddl]&&&& \cel\ar@{-}[ddd] \ar@{-}[ddddl]&{} \\
   10 \ar@{.}[rrrrrrrrrrrrrrrr] &&\cel\ar@{-}[ddddl]\ar@{-}[ddd]  &&&&\cel\ar@{-}[ddddl]\ar@{-}[ddd]
  &&&&\cel\ar@{-}[ddddl] \ar@{-}[ddd] &&&& \cel\ar@{-}[ddddl]\ar@{-}[ddd]  &&{} \\
   9 \ar@{.}[rrrrrrrrrrrrrrrr] &&&&&\cel\ar@{-}[ddd] &&&&\cel\ar@{-}[ddd] &&\cel &&&& &{} \\
   8 \ar@{.}[rrrrrrrrrrrrrrrr] &&&\cel &&&&&&&&&&\cel\ar@{-}[dd] && \cel & {} \\
   7 \ar@{.}[rrrrrrrrrrrrrrrr] &&\cel &&&&\cel &&&&\cel &&&& \cel  &&{} \\
   6 \ar@{.}[rrrrrrrrrrrrrrrr] & \cel &&&&\cel &&&&\cel &&&& \cel &&&{} \\
   && *=0{A(v\eta^2w\eta^2)} &&&& *=0{A(\eta^2 v\eta^2w)} &&&& *=0{A(\eta^2v\eta^2w\eta)}
 &&&& *=0{A(\eta v\eta^2w\eta^2)}&&
 }  
 $$  
 \end{table}

 Again $2$-primary atoms are those with $v,w\in\set{3,6}$ and there are no $3$-primary spaces in this table.
 Moreover, the new $2$-primary atoms are pairwise non-congruent, wherefrom we obtain the following
 result.

 \begin{corollary}\label{tf1-3} 
  The Grothendieck group $K_0(\cwf^5)$ is a free abelian group of rank 85.
 \end{corollary}
 
 We end up with the following statements about the higher dimensional torsion free spaces.

  \begin{proposition}\label{tf1-4} 
  \begin{enumerate}
 \item 	There are infinitely many non-isomorphic (even non-congruent) $2$-primary atoms of dimension $13$.
 Hence the Groth\-en\-dieck group $K_0(\cwf^k)$ is of infinite rank for $k\ge6$.
 \item   If $k\ge11$, the classification problem for the category $\cwf^k$ is wild.
\end{enumerate}
 \end{proposition} 
 \begin{proof} 
  We shall show first that $\pi_{12}^S(A^{11}(\eta^2v)$, or, the same, $\pi_{10}^S(A(\eta^2v)$ equals
 $\mZ/2\+\mZ/2$. We consider the cofibration sequences 
 \begin{align*} 
  \tag{a}& S^8 \str f \Si C \str g A \str h S^9 \str{\Si f} \Si^2C,\\
  \tag{b}& S^6 \larr S^4 \larr C \larr S^7 \larr S^5,
 \end{align*} 
 where $A=A(\eta^2v),\,C=C(\eta^2)$. Note that the map $f$ factors through $S^5$. From the sequence (b)
 we get $\pi_9^S(C)\iso\pi_9^S(S^7)\iso\mZ/2$ and $\pi^S_9(\Si C)\iso\pi^S_8(C)=0$. The second equality
 follows from the fact that the induced map $\pi^S_8(S^7)\to\pi^S_8(S^5)$ is known to be injective \cite{tod}. 
 Since $\pi^S_{10}(S^5)=\pi^S_{10}(S^6)=0$, the sequence (a) gives then an exact sequence
 $$ 
  0 \larr \pi_{10}^S(\Si C)\iso\mZ/2 \larr \pi_{10}^S(A) \larr \pi^S_{10}(S^9)\iso\mZ/2 \larr 0.
 $$ 
 To show that this sequence splits, we have to check that $2\al=0$ for every $\al\in\pi^S_{10}(A)$.
 In any case, $2\al$ factors through $\Si C$, which gives rise to a commutative diagram
 $$ 
   \begin{CD}
    M^{10}(2) @>\phi>> S^{10} @>2>> S^{10} \\
   @V\ga VV	@VV\be V @VV\al V	\\
   S^8  @>>f> \Si C @>>g> A 
 \end{CD}  
 $$ 
 for some $\be,\ga$ (we have used the cofibration sequence for $M(2)$\,). Since $\pi^S_9(S^5)=\pi^S_{10}(S^5)=0$,
 also $\hos(M^{10}(q),S^5)=0$. But the map $\be \phi=f\ga$ factors through $S^5$, so $\be\phi=0$ and $\be=2\si$
 for some $\si\in\pi^S_{10}(\Si C)\iso\mZ/2$. Hence $\be=0$ and $2\al=0$.

 Analogous calculations show that any endomorphism of $A$ acts as a homothety on $\pi_{10}^S(A)$. 
 Since, obviously, $\hos(A,S^{10})=0$, Corollary~\ref{bh-2} shows that the category of spaces arising as
 cones of mappings $mS^{12}\to nA^11(\eta^2v)$ is equivalent to the category of representations of
 the Kronecker quiver $\ti A_1$, or, the same, of diagrams of $\mZ/2$-vector spaces of the shape
 $V_1\rightrightarrows V_2$. But it is well-known that this quiver is of infinite type, i.e. has infinitely many
 non-isomorphic indecomposable representations. Obviously, all corresponding spaces are $2$-primary
 and non-congruent, which proves the claim (1).
 
 The claim (2) follows from the equality $\pi^S_{20}(S^{11})\iso(\mZ/2)^3$. It implies that the category
 of spaces, which are cones of mappings $mS^{20}\to mS^{11}$ is equivalent to that of diagrams
 $\ds\xymatrix{{V_1\,} \ar[r] \ar@<1ex>[r] \ar@<-1ex>[r] & {\ V_2}}$. The latter is well-known to be
 wild.
 \end{proof} 

 Perhaps, the estimate $11$ in the claim (2) of Proposition~\ref{tf1-4} is too big, but at the moment
 I do not know a better one. On the other hand, we can hope that the classification problem for
 $\cws^6$ is still tame.

 \end{document}